# A study of the problem of optimal management of discrete product inventory in the stochastic regeneration model with continuous consumption and random delivery delay

P.V. Shnurkov and N.A. Vakhtanov

The paper considers the optimal control problem of inventory of a discrete product in regeneration scheme with a Poisson flow of customer requirements. In the system deferred demand is allowed, the volume of which is limited by a given value. The control parameter is the level of the stock, at which achievement it is necessary to make an order for replenishment. The indicator of management effectiveness is the average specific profit received in one regeneration period. The optimal control problem is solved on the basis of the statement about the extremum of a fractional-linear integral functional on the set of discrete probability distributions. The significant content of the work is the derivation of analytical representations for the mathematical expectation of the increment of the functional of profit obtained during the regeneration period. The obtained analytical representations enable us to explicitly determine the stationary cost indicator of control efficiency, which was introduced in the first part of the research. Thus, it becomes possible to numerically solve the optimal control problem of inventory in the model under consideration.

## 1 Introduction

Stochastic regeneration models for studying inventory management problems have been the subject of a large number of publications. For example, works [1,2] considered different variants of regeneration models for systems of continuous product inventory management. Work [3] studied a special version of the regeneration model describing a system of continuous product inventory management in which stocks are replenished not immediately but over a period of time called a real replenishment period. The version of the stochastic regeneration model considered in [3] is conceptually connected with the classical deterministic model of inventory management formulated, e.g., in [4].

All the above-mentioned works have considered problems of continuous product inventory management. In models of this type, the set of values of a random process?) describing the stock size in the system represents a certain subset of a set of real numbers. Examples of continuous products in real systems are water, oil and oil products, grain, etc. However, there are a number of important products whose volume is measured in discrete values. Among them are many consumer goods, first of all, consumer electronics, food commodities and some others. Thus, an urgent problem of applied mathematics now is studying mathematical models of discrete product inventory management.

Foreign scientists studying the theory of inventory management have long been using stochastic regeneration models. Models in which the points of regeneration coincide with the points of stock replenishment to a certain fixed level, and the management efficiency indicator represents average unit costs during the regeneration period have been mentioned among the results reported in works [4,5]. However, none of the studies conducted in the last decades have considered the problem of discrete product inventory management in the regeneration scheme.

## 2 General description of operation of the discrete inventory management system under consideration

Let us consider a certain trading system (a warehouse) for temporary storage and consumer delivery of a certain product of the same type whose size is measured in discrete integer values. It is assumed that at the start time the system is full and contains $N$ product units. The product is consumed at random moments when orders (customers) arrive. The arrivals form a Poisson flow with the known parameter $\lambda > 0$. At the point of the next order arrival, one product unit is consumed. Such consumption is immediate.

Stocks in the considered system are replenished at certain time intervals in the following way. The next product lot is ordered when the amount of product at a warehouse reaches level $r$, where $r \leq N$ is a certain integer value that will further act as a solution (management). Note that parameter $r$ is probabilistic in nature. A formal description of the selection procedure for this parameter is given in section 3. Let us call the time period from the product order placement till the arrival of a new lot a period of delivery delay or simply a delay period. The duration of this period is a random variable with a given distribution $H_r(x)$ which, generally speaking, can depend on the control parameter $r$. The stock is replenished as soon as the delay period is over.

During the delivery delay period, product consumption continues according to the described rules. After the physical stock is consumed, this inventory has negative values, which corresponds to the volume of unsatisfied demand or deficit. And every order is registered and satisfied at the next replenishment point. Such orders make up the so-called deferred demand. The deferred demand size in this model cannot exceed a given value $N_0$. If the deferred demand reaches the maximum permissible level $N_0$ before the next replenishment point, all the orders placed after that are lost. It is assumed that stock replenishment is arranged so that the deferred demand is completely satisfied and the product stock is replenished up to the initial level $N$. Further functioning of the system does not depend on the previous events and corresponds to the rules described above. Fig. 1 illustrates a possible evolution of the product stock in the system. Note that the control parameter $r$ can take any integer value within the limits $-N_0 \leq r \leq N$.

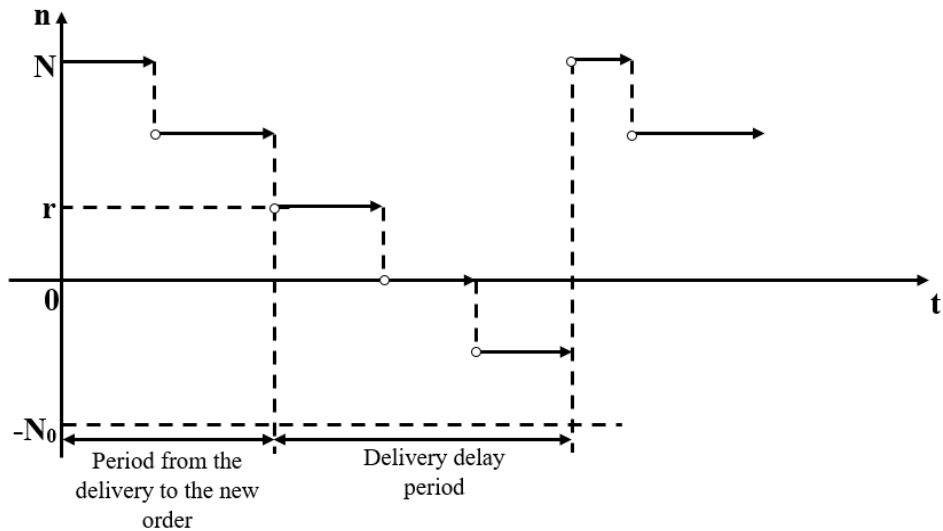

**Fig. 1.** Evolution of product stock in the considered system

## 3  Formal construction of the mathematical model describing how the considered system functions

Let us assume that all the stochastic objects that will be later introduced into the system are determined over a certain probability space $(\Omega, \mathcal{A}, P)$, which represents a formal construction describing a random experiment with a real system. First of all, let us introduce a random process $\xi(t) = \xi(\omega, t), \omega \in \Omega, t \in T = [0, \infty)$ that will act as a mathematical model of the functioning of the considered system. Let us assume that the value of this process at an arbitrary moment of time $t$ represents the product stock in the system, and both the real stock and the deficit, i.e. the value of the deferred demand, are taken into account. Thus, the set of the process $\xi(t)$ values is finite: $X = \{-N_0, -N_0 + 1, \ldots, -1, 0, 1, \ldots, N-1, N\}$.

The moments when the process $\xi(t)$ changes states are the moments of consumption (placement of consumers' orders), if the value of the deferred demand does not exceed the value $N_0$, and the moments of the real stock replenishment. To avoid any ambiguity, let us assume that the process trajectories are continuous on the right.

Let us denote random stock replenishment points by $t_n = t_n(\omega), \omega \in \Omega, n = 0, 1, 2, \ldots, t_0 = 0$. According to our hypothesis about this system's behavior described in the previous section, $\xi(t_n) = N, n = 0, 1, 2, \ldots$ Besides, after every replenishment point, the system evolution continues regardless of the previous events and by the same rules. It means that the process $\xi(t)$ is regenerative [6], and the points $t_n, n = 0, 1, 2, \ldots$ are its regeneration points. The regeneration periods of this process consist of two independent components:

$$\Delta t_n = t_{n+1} - t_n = \Delta_n^{(0)} + \Delta_n^{(1)}, \tag{1}$$

where $\Delta_n^{(0)}$ is a random time from the moment of the previous replenishment till the moment of the next order placement, and $\Delta_n^{(1)}$ is random duration of the delivery delay (the time from the order placement point till the next stock replenishment point).

According to the assumptions made to satisfy the additional condition that the control parameter takes a fixed value $r$, the random variable $\Delta_n^{(0)}$ has Erlang distribution of the $N - r$ order, and the random variable $\Delta_n^{(1)}$ – a given distribution $H_r(x)$.

Let us now describe the procedure of controlling the random process $\xi(t)$.

The control parameter will be the stock level $r$, at which the next product lot is ordered. The control decision is taken at the initial point of every regeneration period, after the next stock replenishment, i.e. at the moments of time $t_n + 0, n = 0, 1, 2, \ldots$

This decision allows determining the level $r$, at which the next product order is placed. Thus, the set of permissible values of control parameter $U$ becomes finite:

$$r \in U = \{N, N-1, \ldots, -N_0\}.$$

Let us define the totality of all possible discrete probability distributions $\Gamma_d$ on set $U$. Every discrete probability distribution belonging to this totality represents a vector $\alpha = (\alpha_N, \alpha_{N-1}, \ldots, \alpha_0, \alpha_{-1}, \ldots, \alpha_{-N_0})$, satisfying the following conditions: $\alpha_i \geq 0, i \in U, \sum_{i \in U} \alpha_i = 1$. Any fixed distribution $\alpha$ defines a probabilistic control strategy, which implies that at every decision-making point, the control value is selected in accordance with this distribution. In other words, the control parameter value $r$ is selected with the probability $\alpha_r, r \in U$. Note that the control strategy is the same for all regeneration periods. The optimization problem consists in finding a control strategy $\alpha^* \in \Gamma_d$ that affords an extremum to a certain stationary cost efficiency indicator.

This indicator depends on the initial cost characteristics of a model. Let us enumerate these characteristics that are assumed to be known:

income from the sales of a product unit $c_0$;

expenses on the storage of a product unit per time unit $c_1$;

expenses on purchasing a product unit $c_2$;

expenses related to the deficit of a product unit per time unit $c_3$;

expenses related to the loss of $i$ clients, $c_4^{(i)}, i = 1, 2, \ldots$

The introduced cost parameters determine all main types of income and expenses appearing in inventory management models.

## 4 Statement of the optimal management problem

First of all, let us introduce an additive cost functional for the considered model. Let $\gamma(t) = \gamma(\omega, t), \omega \in \Omega, t \in T = [0, \infty)$ be random profit gained as a result of functioning of the stochastic system described above over the time interval $[0, t]$, and $\gamma(0) = \gamma_0$ be a specified value. The process $\gamma(t)$ depends on the trajectory of the process $\xi(t)$ and its values are also determined by the initial cost parameters $c_0, c_1, c_2, c_3, c_4^{(i)}, i = 1, 2, \ldots$ A complete description of the trajectories of the process $\gamma(t)$ is complex and not necessary for solving the problem of optimal management in question.

Similar additive cost functionals have been earlier reported in the relevant scientific literature. The general scheme of their construction for controlled Markov and semi-Markov random processes is given in classical works [7, 8]. This work requires studying probabilistic characteristics related to the process $\gamma(t)$ increments over the regeneration periods of the main process $\xi(t)$.

Let us denote by $\gamma_n = \gamma(t_n), n = 0, 1, 2, \ldots$, the values of the process $\gamma(t)$ at the regeneration points assuming that $\gamma(t_n + 0) = \gamma(t_n)$, $n = 0, 1, 2, \ldots$ Let $\Delta\gamma_n = \gamma_{n+1} - \gamma_n$ be the profit increment of the regeneration period $(t_n, t_{n+1}]$, $n = 0, 1, 2, \ldots$ Let us denote by $E_\alpha(\Delta\gamma_n)$ the mathematical expectation of the profit increment over the regeneration period if the strategy of the controlling process $\xi(t)$ is determined by the distribution $\alpha \in \Gamma_d$. In the same way, the value $E_\alpha(\Delta t_n) = E_\alpha(t_{n+1} - t_n)$ represents a mathematical expectation of the regeneration period duration for the control strategy $\alpha \in \Gamma_d$. Let us also introduce notation $E_\alpha \gamma(t)$ for the mathematical expectation of the total profit gained over the time period $[0, t]$ for the control strategy $\alpha \in \Gamma_d$.

In the assumption that the decision over the regeneration period is fixed, i.e. the control parameter assumes the value $r \in U$, let us express the corresponding conditional mathematical expectations as $E_r(\Delta\gamma_n)$ and $E_r(\Delta t_n), n = 0, 1, 2, \ldots$

The change in the additive functional $\gamma(t)$ in one regeneration period $(t_n, t_{n+1}]$ consists of two components:

$$\Delta\gamma_n = \Delta\gamma_n^{(0)} + \gamma_n^{(1)}, \tag{2}$$

where $\Delta\gamma_n^{(0)} = \gamma(t_{n+1} - 0) - \gamma(t_n)$ is the increment of the functional $\gamma(t)$ over the open interval $(t_n, t_{n+1})$, in this model this increment depends on the trajectory of the process $\xi(t)$ over this interval and is determined by the parameters $c_0, c_1, c_3, c_4^{(i)}, i = 1, 2, \ldots$; $\gamma_n^{(1)} = \gamma(t_{n+1}) - \gamma(t_{n+1} - 0)$ is the functional $\gamma(t)$ increment at the time moment $t_{n+1}$, i.e. its instantaneous change at this time moment. In the considered model, the value $\gamma_n^{(1)}$ represents the costs related to the stock replenishment. It means that this variable is negative and its value depends on the stock replenishment size at the point $t_{n+1}$ and is defined by the equality $\gamma_n^{(1)} = -c_2[\xi(t_{n+1}) - \xi(t_{n+1} - 0)], n = 0, 1, 2, ...$

Given the mentioned features of the model, the value of the process $\gamma(t)$ at an arbitrary point of time $t > 0$ also consists of two components:

$$\gamma(t) = \gamma^{(0)}(t) + \gamma^{(1)}(t), \tag{3}$$

where $\gamma^{(0)}(t)$ is the accumulated profit over the interval $[0, t]$, determined by changes in the process $\xi(t)$ over the regeneration periods, excluding the regeneration points themselves if they belong to the interval $[0, t]$; $\gamma^{(1)}(t)$ is the accumulated costs over the interval $[0, t]$ related to the stock replenishment at the regeneration points of the process $\xi(t)$ belonging to the interval $[0, t]$.

Normally, starting from the first fundamental works on the theory of Markov and semi-Markov processes control [7, 8], the control efficiency indicator has been the variable

$$I_\alpha = \lim_{t \to \infty} \frac{E_\alpha \gamma(t)}{t}, \tag{4}$$

which conceptually represents average specific profit in the considered stochastic model defined by control strategy $\alpha$. However, the authors of works [7, 8] studied additive cost functionals of the $\gamma^{(0)}(t)$ type, which did not account for the components appearing at certain random points of time. In order to use the efficiency indicator of type (4), we have to prove the ergodic theorem about the existence and explicit expression of the specified limit for the additive functional of type (3) taking into account the component $\gamma^{(1)}(t)$. Let us prove the corresponding result for the regenerating process $\xi(t)$ and the additive cost functional $\gamma(t)$ determined on it.

**Theorem 1.** Let $\xi(t)$ be a controlled regenerating process with the control strategy determined by the probability distribution $\alpha \in \Gamma_d$, and $\gamma(t)$ be the additive cost functional related to the process $\xi(t)$, whose general structure is defined by relation (3). Let us assume that for any permissible control strategy $\alpha \in \Gamma_d$ the following condition $E_\alpha[\Delta t_n] > 0$ is fulfilled. Then there is a limit

$$I_\alpha = \lim_{t\to\infty} \frac{E_\alpha \gamma(t)}{t} = \frac{E_\alpha[\Delta\gamma_n]}{E_\alpha[\Delta t_n]}, \qquad (5)$$

where $\Delta\gamma_n$ is the complete change in the functional $\gamma(t)$ over the regeneration period $(t_n, t_{n+1}]$, defined by relation (2).

**Proof.** According to the classical renewal theory [9], let us introduce $v(t) = \sup\{n: t_n \leq t\}$. The process $v(t)$ is usually called the counting process for the renewal process formed by the regeneration points $\{t_n, n = 0, 1, 2, ..\}$. Let $H(t) = Ev(t)$ be the function of the regeneration process $\{t_n, n = 0, 1, 2, ..\}$. For a fixed value $t$, the random variable $v(t)$ coincides with the number of renewals over the interval $[0, t]$ (the point $t = 0$ is not considered to be a point of renewal).

Based on the process $\gamma(t)$ additivity, it can be stated that

$$\gamma^{(1)}(t) = \sum_{n=0}^{v(t)-1} \gamma_n^{(1)} = \sum_{n=1}^{v(t)} \gamma_{n-1}^{(1)}. \qquad (6)$$

Note that the properties of the regeneration process $\xi(t)$ and process $\gamma(t)$ related to it suggest that for a random point of time $t > 0$ and any fixed value $n = 1, 2, ...$, the random event $(v(t) \leq n)$ does not depend on the system of events generated by the random variables $\{\gamma_{n+1}^{(1)}, \gamma_{n+2}^{(1)}, ...\}$. Such property in the probability theory is called independence from the future [10]. Indeed, the random variable $v(t)$ and event $(v(t) \leq n)$ depend on random points of time $t_1, t_2, ... t_n, t_{n+1}$ and corresponding regeneration periods $[0, t_1], (t_1, t_2], ..., (t_n, t_{n+1}]$. At the same time, the random variables $\{\gamma_{n+1}^{(1)}, \gamma_{n+2}^{(1)}, ...\}$ are related to the regeneration periods $(t_{n+1}, t_{n+2}], (t_{n+2}, t_{n+3}], ...$

If the property of independence from the future holds, the Wald identity is satisfied [10, ch. 4, § 4]. Based on this result, we derive the following from (6)

$$E_\alpha \gamma^{(1)}(t) = E_\alpha v(t) E_\alpha \gamma_n^{(1)} = H(t) E_\alpha \gamma_n^{(1)}. \qquad (7)$$

From relation (3) and taking into account (7) we get

$$I_\alpha = \lim_{t\to\infty} \frac{E_\alpha \gamma(t)}{t} = \lim_{t\to\infty} \frac{E_\alpha \gamma^{(0)}(t)}{t} + \lim_{t\to\infty} \frac{E_\alpha \gamma^{(1)}(t)}{t} =$$

$$= \lim_{t\to\infty} \frac{E_\alpha \gamma^{(0)}(t)}{t} + E_\alpha \gamma_n^{(1)} \lim_{t\to\infty} \frac{H(t)}{t}. \qquad (8)$$

The first summand on the right side of equality (8) can be determined based on the ergodic theorem for the additive functional. The classical forms of such theorem for semi-Markov models are given in [7, 8]. In this case, the ergodic theorem for the regenerating process $\xi(t)$ brings us to the following result:

$$\lim_{t\to\infty} \frac{E_\alpha \gamma^{(0)}(t)}{t} = \frac{E_\alpha[\Delta \gamma_n^{(0)}]}{E_\alpha[\Delta t_n]}. \tag{9}$$

The second summand on the right side of equality (8) is determined based on the elementary renewal theorem [9]:

$$\lim_{t\to\infty} \frac{H(t)}{t} = \frac{1}{E_\alpha[\Delta t_n]}. \tag{10}$$

Note that we can use this theorem if $E_\alpha[\Delta t_n] > 0$.

Then from (8), taking into account (9) and (10), we derive the following expression for the average specific profit:

$$I_\alpha = \lim_{t\to\infty} \frac{E_\alpha \gamma(t)}{t} = \frac{E_\alpha[\Delta \gamma_n^{(0)}] + E_\alpha \gamma_n^{(1)}}{E_\alpha[\Delta t_n]} = \frac{E_\alpha[\Delta \gamma_n^{(0)} + \gamma_n^{(1)}]}{E_\alpha[\Delta t_n]}. \tag{11}$$

The random variable $\Delta \gamma_n^{(0)} + \gamma_n^{(1)} = \Delta \gamma_n$ is the total fluctuation of the additive cost functional over the regeneration period $(t_n, t_{n+1}]$. The obtained relation (11) coincides with (5). Theorem 1 is proved.

**Note 1.** The statement of theorem 1 holds for a more general variant, in which the set of permissible controls $U$ has an arbitrary structure, and the set of control strategies of the regenerating process $\xi(t)$ coincides with the set of probability distributions defined on $U$.

**Note 2.** The condition $E_\alpha[\Delta t_n] > 0$ can be satisfied in the following way. It follows from relation (1) that $E_\alpha[\Delta t_n] = E_\alpha\left[\Delta_n^{(0)}\right] + E_\alpha\left[\Delta_n^{(1)}\right]$, where $E_\alpha\left[\Delta_n^{(0)}\right] = \sum_{r\in U} \frac{N-r}{\lambda}\alpha_r$, $E_\alpha\left[\Delta_n^{(1)}\right] = \sum_{r\in U} \tau_r \alpha_r$, $\tau_r = E_r\left[\Delta_n^{(1)}\right] = \int_0^\infty x dH_r(x)$.

It is evident that $E_\alpha\left[\Delta_n^{(0)}\right] = 0$ only for the distribution of the following type: ($\alpha_N = 1, \alpha_r = 0, r \in U, r \neq N$). It follows that if $\tau_N > 0$, then $E_\alpha[\Delta t_n] > 0$ for the specified probability distribution ($\alpha_N = 1, \alpha_r = 0, r \in U, r \neq N$), as well as for any other decision-making strategy $\alpha \in \Gamma_d$. But if we require that a stronger natural condition $\tau_r > 0, r \in U$ (the average duration of a delay period cannot be equal to zero at any permissible decision) is satisfied, the following inequality holds: $E_\alpha[\Delta t_n] \geq E_\alpha\left[\Delta_n^{(1)}\right] > 0$.

It follows from the statement of theorem 1 that to find an explicit expression for the efficiency indicator $I_\alpha$ requires determining analytical expressions of mathematical expectations $E_\alpha[\Delta \gamma_n], E_\alpha[\Delta t_n]$.

Let us fix the control parameter $r \in U$ and denote

$$A(r) = E_r\left[\Delta \gamma_n^{(0)} + \gamma_n^{(1)}\right] = E_r[\Delta \gamma_n]; \tag{12}$$

$$B(r) = E_r[\Delta t_n]. \tag{13}$$

According to the property of the mathematical expectation, the following formulae hold:

$$E_\alpha[\Delta \gamma_n] = \sum_{r \in U} E_r[\Delta \gamma_n] \alpha_r = \sum_{r \in U} A(r) \alpha_r, \tag{14}$$

$$E_\alpha[\Delta t_n] = \sum_{r \in U} E_r[\Delta t_n] \alpha_r = \sum_{r \in U} B(r) \alpha_r. \tag{15}$$

Taking into account (14) and (15) we derive from (11) that

$$I_\alpha = \frac{E_\alpha[\Delta \gamma_n]}{E_\alpha[\Delta t_n]} = \frac{\sum_{r \in U} A(r) \alpha_r}{\sum_{r \in U} B(r) \alpha_r}. \tag{16}$$

The problem of optimal control in the considered model can be formulated as the following extremum problem:

$$I_\alpha = \frac{\sum_{r \in U} A(r) \alpha_r}{\sum_{r \in U} B(r) \alpha_r} \to max, \alpha \in \Gamma_d. \tag{17}$$

Functional (16) in its form is a so-called linear-fractional integral-type functional defined on a set of discrete probability distributions $\Gamma_d$. The general approach to solving this problem is described below.

## 5 Extremum problem for the linear-fractional functional defined on a set of discrete probability distributions

Let us define a certain number of discrete sets $U_i = 1, 2, \ldots, n_i, i = 1, 2, \ldots, N, N < \infty$. Every set $U_i$ can include a finite or denumerable number of elements: $n_i \leq \infty, i = 1, 2, \ldots, N$. These sets can be understood as totalities of permissible decisions (controls) made in different states of the considered stochastic model. On each set $U_i$, we define a totality of possible probabilistic distributions of the following type: $\alpha^{(i)} =$

$\left(\alpha_1^{(i)}, \alpha_2^{(i)}, \ldots, \alpha_{n_i}^{(i)}\right), \alpha_s^{(i)} \geq 0, s = 1, 2, \ldots, n_i; \sum_{s=1}^{n_i} \alpha_s^{(i)} = 1, i = 1, 2, \ldots, N$. Let us denote by $\Gamma_d^{(i)}, i = 1, 2, \ldots, N$ the described totality of probability distributions defined on the set $U_i$. Let us consider the Cartesian product of spaces $U = U_1 \times U_2 \times \ldots \times U_N$. Let us define the probability measure on $U$ as a product of probability measures on the spaces $U_1, U_2, \ldots, U_N$, defined by the distributions $\alpha^{(1)}, \alpha^{(2)}, \ldots, \alpha^{(N)}$. Thus, the probability measure on $U$ is defined by a totality of probability distributions $\alpha^{(1)}, \alpha^{(2)}, \ldots, \alpha^{(N)}$. Let us denote by $\Gamma_d$ the set of probability measures on space $U$ defined by the method described above. Further, additional conditions will be applied to the set $\Gamma_d$ when the extremum problem is set. In advance let us introduce the notion of degenerate discrete probability distribution.

**Definition 1.** Let us term the probability distribution $\alpha^{(i)*}(k_i)$ degenerate if $\alpha_{k_i}^{(i)} = 1, \alpha_l^{(i)} = 0, l = 1, 2, \ldots, n_i, l \neq k_i$ and the point $k_i \in U_i$ a point of degenerate discrete distribution $\alpha^{(i)*}(k_i)$. Degenerate distribution is known to correspond to a deterministic variable that assumes value $k_i$.

The set of degenerate probability distributions defined on $U_i$ will be denoted by the symbol $\Gamma_d^{(i)*}, i = 1, 2, \ldots, N$. There is a one-to-one correspondence between the $\Gamma_d^{(i)*}$ and $U_i, i = 1, 2, \ldots, N$ sets. Let us denote the totality of all degenerate probability measures defined on set $U$ by $\Gamma_d^*$. Thus, every degenerate measure from the $\Gamma_d^*$ set is defined by a certain number of degenerate distributions $\alpha^{(1)*}, \alpha^{(2)*}, \ldots, \alpha^{(N)*}$.

Now let two numerical functions be defined on set $U$:

$A(k_1, k_2, \ldots, k_N): U \to R, B(k_1, k_2, \ldots, k_N): U \to R$, где $k_i \in U_i, i = 1, 2, \ldots, N$.

The integral with respect to a certain discrete measure defined on a discrete set is known to transform into a sum. The corresponding multidimensional integral with respect to a measure generated by the product of the initial measures defined on the Cartesian product of discrete spaces is a multidimensional sum. Now, let us define the linear-fractional integral-type functional and the main function.

**Definition 2.** The linear-fractional integral-type functional in discrete representation or linear-fractional integral-type functional defined on the set of series of discrete probability distributions $\Gamma_d$, is the mapping $I(\alpha^{(1)}, \alpha^{(2)}, \ldots, \alpha^{(N)}): \Gamma_d \to R$, defined by the following equality:

$$I(\alpha^{(1)}, \alpha^{(2)}, \ldots, \alpha^{(N)}) = \frac{\sum_{k_1=1}^{n_1} \sum_{k_2=1}^{n_2} \cdots \sum_{k_N=1}^{n_N} A(k_1, k_2, \ldots, k_N) \alpha_{k_1}^{(1)} \alpha_{k_2}^{(2)} \cdots \alpha_{k_N}^{(N)}}{\sum_{k_1=1}^{n_1} \sum_{k_2=1}^{n_2} \cdots \sum_{k_N=1}^{n_N} B(k_1, k_2, \ldots, k_N) \alpha_{k_1}^{(1)} \alpha_{k_2}^{(2)} \cdots \alpha_{k_N}^{(N)}} \quad (18)$$

**Definition 3.** The main function of the linear-fractional integral-type functional $I(\alpha^{(1)}, \alpha^{(2)}, \ldots, \alpha^{(N)})$ is the function $C(k_1, k_2, \ldots, k_N): U \to R$, defined by the equality

$$C(k_1, k_2, \ldots, k_N) = \frac{A(k_1, k_2, \ldots, k_N)}{B(k_1, k_2, \ldots, k_N)}. \quad (19)$$

Let us now formulate the extremum problem for the introduced linear functional $I(\alpha^{(1)}, \alpha^{(2)}, \ldots, \alpha^{(N)})$ of type (18) on the set of series of discrete probability distributions $\Gamma_d$:

$$I(\alpha^{(1)}, \alpha^{(2)}, \ldots, \alpha^{(N)}) \to extr, (\alpha^{(1)}, \alpha^{(2)}, \ldots, \alpha^{(N)}) \in \Gamma_d. \quad (20)$$

Similarly to [11], let us introduce several preliminary conditions related to linear-fractional integral-type functional (18) and corresponding to extremum problem (20).

Condition 1. For any series of probability distributions $(\alpha^{(1)}, \alpha^{(2)}, \ldots, \alpha^{(N)}) \in \Gamma_d$ there are functionals which define the numerator and denominator of the functional $I(\alpha^{(1)}, \alpha^{(2)}, \ldots, \alpha^{(N)})$. In other words, it is supposed that the numerical series on the right side of the expressions

$$I_1(\alpha^{(1)}, \alpha^{(2)}, \ldots, \alpha^{(N)}) = \sum_{k_1=1}^{n_1} \sum_{k_2=1}^{n_2} \cdots \sum_{k_N=1}^{n_N} A(k_1, k_2, \ldots, k_N) \alpha_{k_1}^{(1)} \alpha_{k_2}^{(2)} \cdots \alpha_{k_N}^{(N)},$$

$$I_2(\alpha^{(1)}, \alpha^{(2)}, \ldots, \alpha^{(N)}) = \sum_{k_1=1}^{n_1} \sum_{k_2=1}^{n_2} \cdots \sum_{k_N=1}^{n_N} B(k_1, k_2, \ldots, k_N) \alpha_{k_1}^{(1)} \alpha_{k_2}^{(2)} \cdots \alpha_{k_N}^{(N)}$$

are convergent.

Condition 2. The functional $I_2(\alpha^{(1)}, \alpha^{(2)}, \ldots, \alpha^{(N)})$ does not vanish for all the series of discrete distributions $(\alpha^{(1)}, \alpha^{(2)}, \ldots, \alpha^{(N)})$ from the set $\Gamma_d$:

$$I_2(\alpha^{(1)}, \alpha^{(2)}, \ldots, \alpha^{(N)}) = \sum_{k_1=1}^{n_1} \sum_{k_2=1}^{n_2} \cdots \sum_{k_N=1}^{n_N} B(k_1, k_2, \ldots, k_N) \alpha_{k_1}^{(1)} \alpha_{k_2}^{(2)} \cdots \alpha_{k_N}^{(N)} \neq 0.$$

Condition 3. The $\Gamma_d$ set contains the whole set of the series of degenerate probability distributions $\Gamma_d^*$, i.e. $\Gamma_d$: $\Gamma_d^* \subset \Gamma_d$.

**Note 3.** It directly follows from conditions 2 and 3 that the function $B(k_1, k_2, .., k_N)$ must have a strictly constant sign for any $(k_1, k_2, .., k_N) \in U$. However, if this condition for the function $B(k_1, k_2, .., k_N)$ holds, condition 2 will hold automatically. Therefore, in the statement of extremum problem (19) given below, we will require satisfaction of the condition that the function $B(k_1, k_2, .., k_N)$ has strictly the same sign for all $(k_1, k_2, .., k_N) \in U$.

**Definition 4.** The totality of the series of discrete probability distributions $\Gamma_d$ is called acceptable in the extremum problem (20) set if Conditions 1 and 3 of the above system of preliminary conditions hold.

Let us now formulate the main statement of this section about solving extremum problem (20).

**Theorem 2.** Let the set of the series of discrete probability distributions $\Gamma_d$ be permissible with respect to the extremum problem (20). Suppose the function $B(k_1, k_2, .., k_N)$ in the linear-fractional integral-type discrete functional $I(\alpha^{(1)}, \alpha^{(2)}, .., \alpha^{(N)})$ has a strictly constant sign (strictly positive or strictly negative) for any values of the arguments $(k_1, k_2, .., k_N) \in U$. Also suppose that the main function $C(k_1, k_2, .., k_N) = \frac{A(k_1, k_2, .., k_N)}{B(k_1, k_2, .., k_N)}$ of the linear-fractional integral-type discrete functional reaches a global extremum (maximum or minimum) on the set $U$ at a certain fixed point $(k_1^*, k_2^*, .., k_N^*)$. Then there is a solution to the set extremum problem and this solution is achieved on a series of degenerate probability distributions $(\alpha^{(1)*}, \alpha^{(2)*}, .., \alpha^{(N)*})$ concentrated at the points $k_1^*, k_2^*, .., k_N^*$, respectively, with the following relations satisfied:

$$\max_{(\alpha^{(1)}, \alpha^{(2)}, .., \alpha^{(N)}) \in \Gamma_d} I(\alpha^{(1)}, \alpha^{(2)}, .., \alpha^{(N)}) =$$

$$= \max_{(\alpha^{(1)*}, \alpha^{(2)*}, .., \alpha^{(N)*}) \in \Gamma_d^*} I(\alpha^{(1)*}, \alpha^{(2)*}, .., \alpha^{(N)*}) =$$

$$= \max_{(k_1, k_2, .., k_N) \in U} \frac{A(k_1, k_2, .., k_N)}{B(k_1, k_2, .., k_N)} = \frac{A(k_1^*, k_2^*, .., k_N^*)}{B(k_1^*, k_2^*, .., k_N^*)},$$

if the main function reaches a global maximum at the point $(k_1^*, k_2^*, .., k_N^*)$;

$$\min_{(\alpha^{(1)}, \alpha^{(2)}, ..., \alpha^{(N)}) \in \Gamma_d} I(\alpha^{(1)}, \alpha^{(2)}, .., \alpha^{(N)}) =$$

$$= \min_{(\alpha^{(1)*}, \alpha^{(2)*}, ..., \alpha^{(N)*}) \in \Gamma_d^*} I(\alpha^{(1)*}, \alpha^{(2)*}, .., \alpha^{(N)*}) =$$

$$= \min_{(k_1, k_2, ..., k_N) \in U} \frac{A(k_1, k_2, .., k_N)}{B(k_1, k_2, .., k_N)} = \frac{A(k_1^*, k_2^*, .., k_N^*)}{B(k_1^*, k_2^*, .., k_N^*)},$$

if the main function reaches a global minimum at the point $(k_1^*, k_2^*, .., k_N^*)$.

The statement of theorem 2 follows from a more general statement about the unconditional extremum of the linear-fractional integral-type functional defined on a set of arbitrary probability measures [11].

## 6 On solving the problem of optimal management in the framework of the developed stochastic model

The theoretical result related to solving the problem of unconditional extremum (20) for the linear-fractional integral-type functional defined on a set of discrete probability distributions formulated in the form of theorem 2 can be used for solving problems of optimal management (17). According to this theorem, if the main function of functional (16) defined by the formula $C(r) = \frac{A(r)}{B(r)}$, where the functions $A(r)$ and $B(r)$ are defined by relations (12) and (13), reaches a global maximum at a certain fixed point $r^* \in U = \{N, N-1, .., 0, -1, .., -N_0\}$, there is a solution to the initial problem (17) and this solution is achieved on the degenerate distribution $\alpha^*$, concentrated at point $r^*$. Note that in the considered problem, the set of feasible solutions $U$ is finite, and the global maximum of the function $C(r)$ is reached. Thus, in order to find a complete solution to the problem set, we have to find explicit analytical expressions for the functions $A(r)$ and $B(r)$. The corresponding results are given in section 8 of the present study.

To present these results, we use the notation $I_r = C(r) = \frac{A(r)}{B(r)}$. In essence, this variable is the value of the linear-fractional integral-type functional on the degenerate discrete distribution concentrated at the point $r \in U$. The solution to the set problem of optimal management (17) will be defined by the point of the global maximum of the function $I_r, r \in U$.

# 7 Auxiliary theoretical results

First of all, let us generalize the known formula for a mathematical expectation of the non-negative random variable.

Let $(\Omega, \mathcal{A}, P)$ be the initial probability space. Suppose we define a non-negative random variable $\eta = \eta(\omega), \omega \in \Omega$ and an arbitrary event $A \in \mathcal{A}$ on this space. Denote by $F(x) = P(\eta < x), 0 < x < \infty$ the function of distribution of the random variable $\eta$.

It is known that B.V.Gnedenko [12] obtained the following result of representing the mathematical expectation of the non-negative random variable $\eta$:

$$E\eta = \int_0^\infty x dF(x) = \int_0^\infty P(\eta \geq x) dx = \int_0^\infty (1 - F(x)) dx. \tag{21}$$

Let us now consider the joint distribution $F(x; A) = P(\eta < x; A)$ and the corresponding mathematical expectation

$$E(\eta; A) = \int_0^\infty x dF(x; A).$$

The following statement of representing a mathematical expectation by joint distribution of the non-negative random variable η and event $A$ was proved by P.V. Shnurkov:

$$E(\eta; A) = \int_0^\infty P(\eta \geq x; A) dx. \tag{22}$$

Relation (22) is a generalization of formula (21) which directly follows from (22) if $A = \Omega$.

General formula (22) will be the basis for the proof of the statement which in this work will be referred to as theorem 3. This statement, in turn, will be used to derive explicit expressions for different auxiliary probability characteristics related to the target management indicator.

Let us now formulate and prove the mentioned statement supposing that all the stochastic objects introduced are defined on the initial probability space $(\Omega, \mathcal{A}, P)$ that describes the random experiment conducted.

**Theorem 3.** Let us consider a Poisson flow of homogeneous events with a parameter $\lambda$. Denote by $\xi_1, \ldots, \xi_s$ the intervals between the consecutive events of this flow, by $\varsigma_s = \sum_{i=1}^{s} \xi_i$ – the point when the *s*-th event occurs, $s = 1, 2, \ldots$

Let $h_r$ be a positive random variable with a given distribution function $F_r(z)$ that depends on a certain real parameter $r$. Let us interpret $h_r$ as a random time interval. Suppose the random variable $h_r$ at any fixed $r$ does not depend on the Poisson flow considered and, in particular, on the random variables $\xi_1, \ldots, \xi_s$. Denote by $A_s$ a random event consisting in the occurrence of *s* events of the Poisson flow over the time interval $h_r$ and by $\theta_{r,s} = h_r - \varsigma_s$ a random variable representing the residual length of the time interval $h_r$ after the occurrence of the *s*-th event of the flow. Then the mathematical expectation of the random variable $\theta_{r,s}$ defined by the joint distribution with event $A_s$ can be expressed by the following formula:

$$\tau_{r,s} = E[\theta_{r,s}; A_s] = \int_0^\infty \frac{\lambda^s}{s!} \int_x^\infty (z-x)^s e^{-\lambda z} dF_r(z) dx. \qquad (23)$$

**Proof.** Note that the random variable $\theta_{r,s}$ is defined on a set of elementary outcomes $\omega \in \Omega$ corresponding to event $A_s$. Then we can use formula (22) given above for the mathematical expectation by the joint distribution of the random variable $\theta_{r,s}$ and event $A_s$:

$$\tau_{r,s} = E[\theta_{r,s}; A_s] = \int_0^\infty P(\theta_{r,s} > x, A_s) dx. \qquad (24)$$

At the same time, event $A_s$ can be expressed as follows:

$$A_s = (\varsigma_s = \sum_{i=1}^s \xi_i < h_r, \varsigma_s + \xi_{s+1} > h_r).$$

Let us consider the joint distribution of the random variable $\theta_{r,s}$ and event $A_s$. For any $x \geq 0$ we get

$$P(\theta_{r,s} > x; A_s) = P(\theta_{r,s} > x, \varsigma_s < h_r, \varsigma_s + \xi_{s+1} > h_r) =$$
$$= P(h_r - \varsigma_s > x, \varsigma_s < h_r, \varsigma_s + \xi_{s+1} > h_r) = P(h_r - \varsigma_s > x, \xi_{s+1} > h_r - \varsigma_s). \quad (25)$$

Due to the suppositions introduced, the random variable $\varsigma_s$ has the Erlang distribution of order $s \geq 1$ with parameter $\lambda$. For this variable there are explicit expressions of the distribution function $G_s(y)$ and its total differential:

$$G_s(y) = P(\varsigma_s < y) = 1 - \sum_{i=0}^{s-1} \frac{(\lambda y)^i}{i!} e^{-\lambda y}, s \geq 1;$$

$$dG_s(y) = \lambda \frac{(\lambda y)^{s-1}}{(s-1)!} e^{-\lambda y} dy, s \geq 1.$$

Note that under the conditions of the theorem that is being proved, the random variables $h_r, \varsigma_s, \xi_{s+1}$ are independent. Fix the values of the random variables $h_r = z$ and $\varsigma_s = y$. Under the specified condition, the event $(h_r - \varsigma_s > x, \xi_{s+1} > h_r - \varsigma_s)$ occurs if the relation $(h_r - \varsigma_s > x, \xi_{s+1} > h_r - \varsigma_s)$ holds, i.e. $0 < y < z - x$. At the same time, the relation $z > x$ must hold as otherwise the event $(h_r - \varsigma_s > x)$ is impossible, i.e. its probability is equal to zero.

It follows from these notes that the joint probability $P(h_r - \varsigma_s > x, \xi_{s+1} > h_r - \varsigma_s)$ can be expressed as

$$P(h_r - \varsigma_s > x, \xi_{s+1} > h_r - \varsigma_s) = \int_x^\infty \int_0^{z-x} e^{-\lambda(z-y)} d\, G_s(y) dF_r(z) =$$

$$= \int_x^\infty \int_0^{z-x} e^{-\lambda(z-y)} \lambda \frac{(\lambda y)^{s-1}}{(s-1)!} e^{-\lambda y} dy dF_r(z) = \int_x^\infty \int_0^{z-x} \lambda \frac{(\lambda y)^{s-1}}{(s-1)!} e^{-\lambda z} dy dF_r(z) =$$

$$= \int_x^\infty \frac{\lambda^s}{(s-1)!} e^{-\lambda z} \int_0^{z-x} y^{s-1} dy dF_r(z) = \int_x^\infty \frac{\lambda^s}{(s-1)!} e^{-\lambda z} \frac{(z-x)^s}{s} dF_r(z) =$$

$$= \int_x^\infty \frac{\lambda^s}{s!} (z-x)^s e^{-\lambda z} dF_r(z). \tag{26}$$

From (25) and (26) we get an explicit expression for the joint distribution $P(\theta_{r,s} > x; A_s), x > 0$. And the only other thing to be done is to apply formula (24):

$$\tau_{r,s} = E[\theta_{r,s}; A_s] = \int_0^\infty P(\theta_{r,s} > x; A_s) dx =$$

$$= \int_0^\infty P(h_r - \varsigma_s > x, \xi_{s+1} > h_r - \varsigma_s) dx =$$

$$= \int_0^\infty \int_x^\infty \frac{\lambda^s (z-x)^s}{s!} e^{-\lambda z} dF_r(z) dx = \int_0^\infty \frac{\lambda^s}{s!} \int_x^\infty (z-x)^s e^{-\lambda z} dF_r(z) dx.$$

The obtained relation coincides with formula (23).

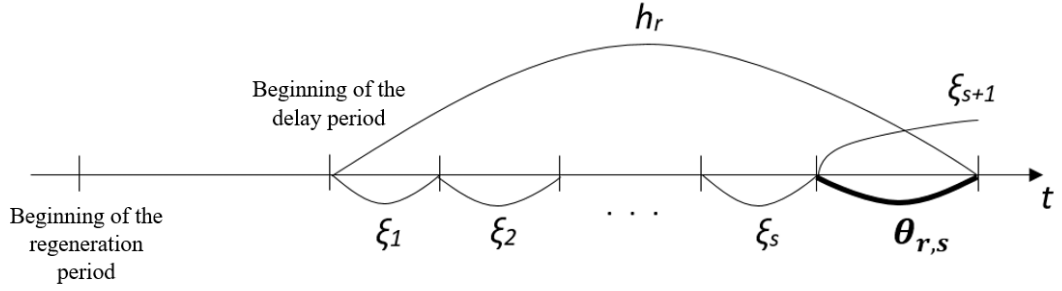

Fig. 2 Graphic illustration of the statement of theorem 3

**Note 1.** In our further analysis of the inventory management model, the Poisson flow of events will be represented by a flow of points of product purchase order arrival. The random variable $h_r$ will play the role of the delivery delay duration with its distribution depending on the control parameter $r$. The random variable $\theta_{r,s}$ will represent the residual duration of the delay period after the arrival of the last order (with a fixed number $s$) in this regeneration period starting from the beginning of the delay period. These events in the process of the evolution of the considered system are illustrated in Fig. 2.

**Note 2.** According to the suppositions made for the model in question, the random variable $h_r$ has a predetermined distribution $H_r(z)$. Thus, further, when using theorem 1, we assume that $F_r(z) = H_r(z), z \geq 0$.

## 8  Analytical expressions for the cost performance of the model under study

Let us move on to proving the statements of specific expressions of probability characteristics related to the stationary cost indicator of management efficiency. First of all, let us introduce a system $\{A_s, s = 0, 1, 2, \dots\}$, where $A_s \in \mathcal{A}$ is an event consisting in the fact that over the delay period the considered system received $s$ orders for the product (clients) exactly. From the Poisson flow properties and formula of total probability it follows that

$$P(A_s) = \int_0^\infty \frac{(\lambda y)^s}{s!} e^{-\lambda y} \, dH_r(y), \qquad s = 0, 1, 2 \ldots$$

Consider the random variable $\Delta\gamma_n$, representing the increment of the additive functional of profit on another regeneration period $(t_n, t_{n+1}]$. Since the system $\{A_s, s = 0, 1, 2, \ldots\}$ forms a complete group of disjoint events, according to the mathematical expectation property it can be written that

$$E(\Delta\gamma_n) = \sum_{s=0}^\infty E(\Delta\gamma_n; A_s).$$

Let us now consider the corresponding mathematical expectations for a fixed value $r$ of the control parameter of this regeneration period. Denote these mathematical expectations $E_r(\Delta\gamma_n), E_r(\Delta\gamma_n; A_s)$. Applying the above-mentioned mathematical expectation property, we get the following initial relation:

$$E_r(\Delta\gamma_n) = \sum_{s=0}^\infty E_r(\Delta\gamma_n; A_s). \tag{27}$$

Let us introduce two random variables related to the regeneration period $(t_n, t_{n+1}]$: $\Delta\gamma_n^{(+)}$ – income increment (in short – income) and $\Delta\gamma_n^{(-)}$ – cost increment (in short – costs). Then the profit increment of a period is represented as $\Delta\gamma_n = \Delta\gamma_n^{(+)} - \Delta\gamma_n^{(-)}$.

Let us also introduce notations $E_r\left(\Delta\gamma_n^{(+)}; A_s\right)$ and $E_r\left(\Delta\gamma_n^{(-)}; A_s\right)$ for conditional mathematical expectations of income and costs over the regeneration period $(t_n, t_{n+1}]$, defined by the joint distribution with event $A_s$ if the solution $r$ is accepted on this period. Since the variable $\Delta\gamma_n$ is additive, the following relation holds

$$E_r(\Delta\gamma_n; A_s) = E_r\left(\Delta\gamma_n^{(+)}; A_s\right) - E_r\left(\Delta\gamma_n^{(-)}; A_s\right). \tag{28}$$

The conditional mathematical expectations $E_r(\Delta\gamma_n)$ are necessary for obtaining an explicit analytical expression of the management efficiency indicator $I_r$.

Note that the formation structure of the random variables $\Delta\gamma_n^{(+)}$ and $\Delta\gamma_n^{(-)}$ strongly depends on the control parameter $r$ and number of orders that arrived over a delivery delay period, i.e. on the value $s$. Bearing this in mind, let us partition the possible values of parameter $r$ into the following disjoint subsets and define an explicit expression for $E_r(\Delta\gamma_n)$ and $I_r$ on each of them:

- $r = N$;
- $1 \leq r < N$;
- $r = 0$;
- $-N_0 < r \leq -1$;
- $r = -N_0$.

Further, in each of the specified values of parameter $r$ it is necessary to find explicit expressions for the mathematical expectations $E_r(\Delta\gamma_n; A_s)$ for all the values $s = 0, 1, 2, \ldots$ After that, we can apply formula (27).

Let us now consider different variants of the values of parameters $r, s$.

Before the analysis, let us make one important methodological note. For ease of reading, let us add commentaries to certain long formulae in the form of special explanatory texts placed under the corresponding parts of this formula. Such texts will be typed in italics.

## 8.1 The variant when r = N

### 8.1.1 Consider the case when s = 0.

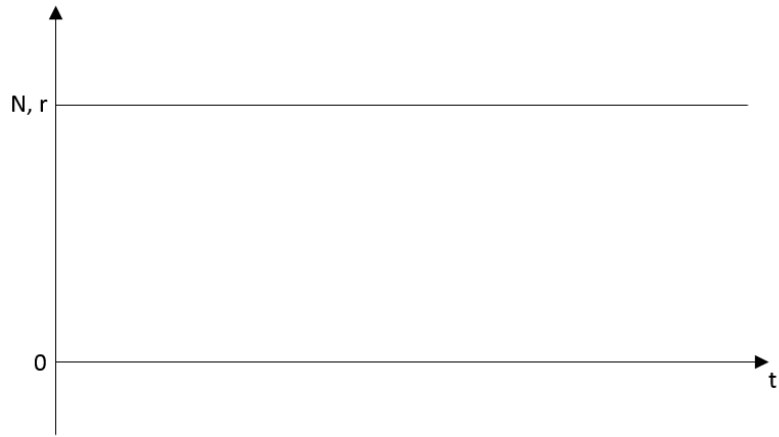

Fig. 3. Graphic illustration for the case 8.1.1

In this variant, an event is realized consisting in the fact that during the regeneration period not a single demand for the product was received. Hence, the mathematical expectation of income by the joint distribution with event $A_0$ is written by the formula

$$E_r\left(\Delta\gamma_n^{(+)}; A_0\right) = 0. \tag{29}$$

The zero value of income is explained by the fact that during the period of regeneration under consideration there were no sales of goods, which is the only source of income for the trading system in question.

The costs in this case coincide with the costs of storing all the goods for the entire regeneration period.

$$E_r(\triangle\gamma_n^{(-)};A_0) = E_r(\triangle\gamma_n^{(-)} \mid A_0)P(A_0) = \underbrace{c_1 N \int_0^\infty (1-H_r(y))dy}_{product\ storage\ costs} \underbrace{\int_0^\infty e^{-\lambda y}dH_r(y)}_{P(A_0)}. \quad (30)$$

Based on formula (28), taking into account (29), (30), we obtain a representation for the mathematical expectation of profit on the regeneration period by the joint distribution with event $A_0$.

$$E_r(\triangle\gamma_n; A_0) = -c_1 N \int_0^\infty (1-H_r(y))dy \int_0^\infty e^{-\lambda y}dH_r(y). \quad (31)$$

Naturally, for this case, when sales on the regeneration period do not occur, the mathematical expectation is negative.

Note that when considering the following options, a similar scheme will be used: expressions for mathematical expectations $E_r(\triangle\gamma_n^{(+)};A_s)$, $E_r(\triangle\gamma_n^{(-)};A_s)$ are determined, and then formula (28) is applied.

8.1.2 Consider the case when $1 \leq s < N$.

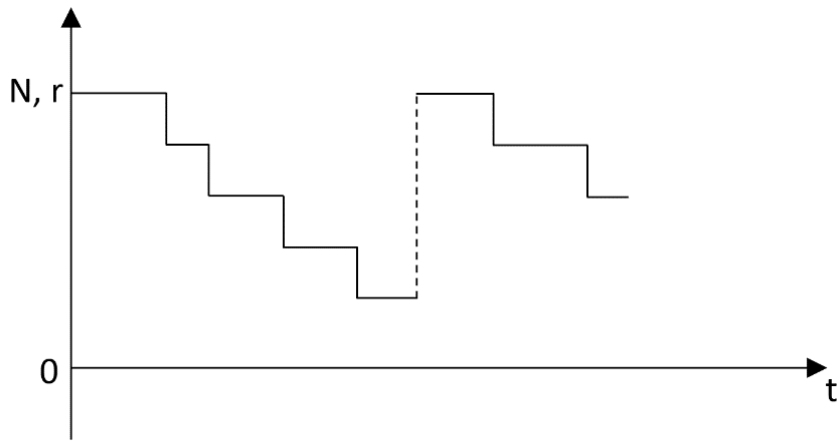

Fig. 4. Graphic illustration for the case 8.1.2

The conditional mathematical expectation of income during the regeneration period, provided that the event $A_s$ occurs, is expressed by the formula

$$E_r(\triangle\gamma_n^{(+)} | A_s) = c_0 s.$$

It follows directly from here that the mathematical expectation of income during the regeneration period by the joint distribution with the event $A_s$ has the form

$$E_r(\triangle\gamma_n^{(+)}; A_s) = E_r(\triangle\gamma_n^{(+)} | A_s)P(A_s) = c_0 s P(A_s) =$$

$$= c_0 s \int_0^\infty \frac{(\lambda y)^s}{s!} e^{-\lambda y} dH_r(y).$$

Let's move on to the analysis of expenses during the regeneration period.

Note that the costs during the regeneration period consist of the costs of storing the product and the costs associated with the acquisition of goods to replenish the stock. Moreover, due to the qualitative difference in the forms of the analytical presentation, it becomes necessary to divide the costs of storing the product during the regeneration period into two components:

1) the part formed on the time interval from the moment the regeneration period begins to the moment the last demand arrives in this period;

2) the part formed on the time interval from the moment of receipt of the last requirement until the end of the regeneration period.

Since such a structure of storage costs and deficit-related penalties will be characteristic of many subsequent options, we introduce a special denotation $E_r^{(1*)}(\triangle\gamma_n^{(-)} | A_s)$ for the conditional mathematical expectation of storage costs arising from the time the last request arrives until the end of the regeneration period, provided that the event $A_s$ occurs. The denotation $E_r^{(1*)}(\triangle\gamma_n^{(-)}; A_s)$ will be used for the corresponding mathematical expectation by the joint distribution with the event $A_s$. From here follows the general formula for the conditional mathematical expectation of expenses during the regeneration period, provided that the event $A_s$ occurs

$$E_r(\triangle\gamma_n^{(-)} | A_s) = \underbrace{\sum_{i=0}^{s-1} \frac{c_1(N-i)}{\lambda}}_{\substack{\text{product storage} \\ \text{costs before} \\ \text{the last order}}} + \underbrace{E_r^{(1*)}(\triangle\gamma_n^{(-)} | A_s)}_{\substack{\text{product storage} \\ \text{costs after} \\ \text{the last order}}} + \underbrace{c_2 s.}_{\substack{\text{stock} \\ \text{replenishment} \\ \text{costs}}}$$

Then, for the corresponding mathematical expectation of the costs by the joint distribution with the event $A_s$, we have

$$E_r(\triangle\gamma_n^{(-)};A_s) = E_r(\triangle\gamma_n^{(-)} | A_s)P(A_s) = \sum_{i=0}^{s-1}\frac{c_1(N-i)}{\lambda}P(A_s) + E_r^{(1*)}(\triangle\gamma_n^{(-)};A_s) + c_2 sP(A_s).$$

According to the formula of the sum of arithmetic progression

$$\sum_{i=0}^{s-1}\frac{c_1(N-i)}{\lambda} = \frac{c_1 s(2N-s+1)}{2\lambda}.$$

To find the mathematical expectation by the joint distribution $E_r^{(1*)}(\triangle\gamma_n^{(-)};A_s)$, we use the statement of Theorem 3 (formula (23)):

$$E_r^{(1*)}(\triangle\gamma_n^{(-)};A_s) = c_1(N-s)\tau_{r,s} = c_1(N-s)\int_0^\infty \frac{\lambda^s}{s!}\int_x^\infty (z-x)^s e^{-\lambda z}dH_r(z)dx.$$

Taking into account the comments made, the mathematical expectation of expenses during the regeneration period by the joint distribution with the event $A_s$ will take the form

$$E_r(\triangle\gamma_n^{(-)};A_s) = \frac{c_1 s(2N-s+1)}{2\lambda}\int_0^\infty \frac{(\lambda y)^s}{s!}e^{-\lambda y}dH_r(y) +$$

$$+c_1(N-s)\int_0^\infty \frac{\lambda^s}{s!}\int_x^\infty (z-x)^s e^{-\lambda z}dH_r(z)dx + c_2 s\int_0^\infty \frac{(\lambda y)^s}{s!}e^{-\lambda y}dH_r(y).$$

Now we can obtain an explicit representation for the mathematical expectation of profit during the regeneration period by the joint distribution with the event $A_s$ for the considered variant of the parameter r, s values

$$E_r(\triangle\gamma_n;A_s) = E_r(\triangle\gamma_n^{(+)};A_s) - E_r(\triangle\gamma_n^{(-)};A_s) =$$

$$c_0 s\int_0^\infty \frac{(\lambda y)^s}{s!}e^{-\lambda y}dH_r(y) - \frac{c_1 s(2N-s+1)}{2\lambda}\int_0^\infty \frac{(\lambda y)^s}{s!}e^{-\lambda y}dH_r(y) -$$

$$-c_1(N-s)\int_0^\infty \frac{\lambda^s}{s!}\int_x^\infty (z-x)^s e^{-\lambda z}dH_r(z)dx - c_2 s\int_0^\infty \frac{(\lambda y)^s}{s!}e^{-\lambda y}dH_r(y) = (c_0 s -$$

$$-\frac{c_1 s(2N-s+1)}{2\lambda} - c_2 s)\int_0^\infty \frac{(\lambda y)^s}{s!}e^{-\lambda y}dH_r(y) - c_1(N-s)\int_0^\infty \frac{\lambda^s}{s!}\int_x^\infty (z-x)^s e^{-\lambda z}dH_r(z)dx. \quad (32)$$

8.1.3 Consider the case when s = N.

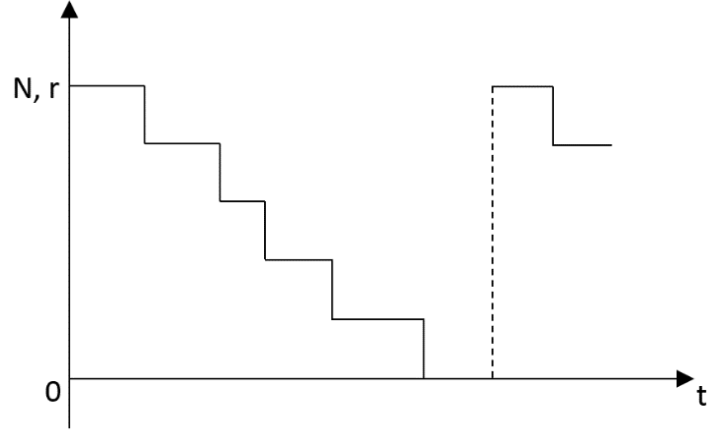

Fig. 5. Graphic illustration for the case 8.1.3

In the case when s = N, the mathematical expectation of income during the regeneration period by the joint distribution with the event $A_N$ is associated with the selling of N units of goods

$$E_r(\triangle\gamma_n^{(+)}; A_N) = E_r(\triangle\gamma_n^{(+)} | A_N)P(A_N) = c_0 N \int_0^\infty \frac{(\lambda y)^N}{N!} e^{-\lambda y} dH_r(y),$$

The corresponding costs are the sum of the storage costs for each unit of the product and the one-time costs of purchasing N units of goods.

$$E_r(\triangle\gamma_n^{(-)}; A_N) = E_r(\triangle\gamma_n^{(-)} | A_N)P(A_N) = (\underbrace{\sum_{i=0}^{N-1}\frac{c_1(N-i)}{\lambda}}_{\text{product storage costs}} + \underbrace{c_2 N}_{\substack{\text{stock} \\ \text{replenishment} \\ \text{costs}}})P(A_N).$$

According to the formula of the sum of arithmetic progression

$$\sum_{i=0}^{N-1}\frac{c_1(N-i)}{\lambda} = \frac{c_1 N(N+1)}{2\lambda}.$$

Thus,

$$E_r(\triangle\gamma_n^{(-)}; A_N) = (\frac{c_1 N(N+1)}{2\lambda} + c_2 N)\int_0^\infty \frac{(\lambda y)^N}{N!} e^{-\lambda y} dH_r(y).$$

Now we can write an explicit representation for the mathematical expectation of profit during the regeneration period by its joint distribution with the event $A_N$

$$E_r(\triangle\gamma_n; A_N) = (c_0 N - \frac{c_1 N(N+1)}{2\lambda} - c_2 N)\int_0^\infty \frac{(\lambda y)^N}{N!} e^{-\lambda y} dH_r(y). \qquad (33)$$

8.1.4 Consider the case when $N < s < N + N_0$.

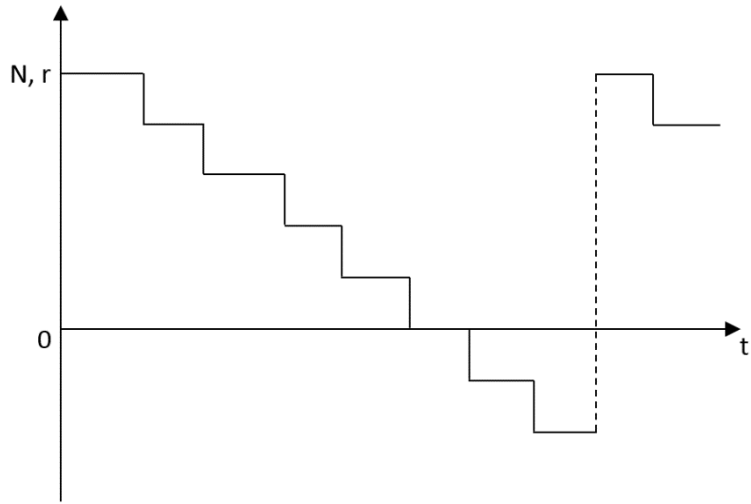

Fig. 6. Graphic illustration for the case 8.1.4

As before, we characterize the structure of the conditional mathematical expectation of expenses provided that the event $A_s$ occurs. They consist of the product storage costs, the stock replenishment costs and the costs (penalties) associated with the deficit of the product. In this case, it is necessary to separate the penalties associated with the deficit and arising in the interval from the moment of the complete exhaustion of the stock until the last request is received on this regeneration period, and the penalties associated with the deficit and arising in the time interval from the moment the last demand is received until the end of the regeneration period. By analogy with the previous case, we will use the notation $E_r^{(3*)}(\triangle \gamma_n^{(-)} | A_s)$ for the conditional mathematical expectation of the second component of penalties.

Thus, the following representation holds for the conditional mathematical expectation of costs

$$E_r(\triangle \gamma_n^{(-)} | A_s) = \underbrace{\sum_{i=0}^{N-1} \frac{c_1(N-i)}{\lambda}}_{\substack{\text{product storage} \\ \text{costs}}} + \underbrace{c_2 s}_{\substack{\text{stock} \\ \text{replenishment} \\ \text{costs}}} + \underbrace{\sum_{i=1}^{s-N-1} \frac{c_3 i}{\lambda}}_{\substack{\text{costs related to the deficit before} \\ \text{the arrival of the last order}}} + \underbrace{E_r^{(3*)}(\triangle \gamma_n^{(-)} | A_s)}_{\substack{\text{costs related to the deficit after} \\ \text{the arrival of the last order}}}.$$

We pass from the conditional mathematical expectation to the corresponding mathematical expectation by the joint distribution

$$E_r(\triangle \gamma_n^{(-)}; A_s) = E_r(\triangle \gamma_n^{(-)} | A_s)P(A_s) =$$

$$= \sum_{i=0}^{N-1} \frac{c_1(N-i)}{\lambda} P(A_s) + c_2 s P(A_s) + P(A_s) \sum_{i=1}^{s-N-1} \frac{c_3 i}{\lambda} + E_r^{(3*)}(\triangle \gamma_n^{(-)}; A_s).$$

To find the mathematical expectation by the joint distribution $E_r^{(3*)}(\triangle \gamma_n^{(-)}; A_s)$, we use the statement of Theorem 3 (formula (23)):

$$E_r^{(3*)}(\triangle \gamma_n^{(-)}; A_s) = c_3(s-N)\tau_{r,s} = c_3(s-N)\int_0^\infty \frac{\lambda^s}{s!} \int_x^\infty (z-x)^s e^{-\lambda z} dH_r(z) dx.$$

According to the formula of the sum of arithmetic progression

$$\sum_{i=0}^{N-1} \frac{c_1(N-i)}{\lambda} = \frac{c_1 N(N+1)}{2\lambda} \quad u \quad \sum_{i=1}^{s-N-1} \frac{c_3 i}{\lambda} = \frac{c_3(s-N)(s-N-1)}{2\lambda}$$

Taking into account the above intermediate results, we obtain

$$E_r(\triangle \gamma_n^{(-)}; A_s) = (\frac{c_1 N(N+1)}{2\lambda} + c_2 s + \frac{c_3(s-N)(s-N-1)}{2\lambda}) \int_0^\infty \frac{(\lambda y)^s}{s!} e^{-\lambda y} dH_r(y) +$$

$$+ c_3(s-N)\int_0^\infty \frac{\lambda^s}{s!} \int_x^\infty (z-x)^s e^{-\lambda z} dH_r(z) dx.$$

The mathematical expectation of revenues by the joint distribution with the event $A_s$ in this case is given by the relation

$$E_r(\triangle \gamma_n^{(+)}; A_s) = E_r(\triangle \gamma_n^{(+)} | A_s)P(A_s) = c_0 s \int_0^\infty \frac{(\lambda y)^s}{s!} e^{-\lambda y} dH_r(y).$$

Thus, the corresponding mathematical expectation of profit for the analyzed variant of the parameter r, s values has the form

$$E_r(\triangle \gamma_n; A_s) = E_r(\triangle \gamma_n^{(+)}; A_s) - E_r(\triangle \gamma_n^{(-)}; A_s) =$$

$$= (c_0 s - \frac{c_1 N(N+1)}{2\lambda} - c_2 s - \frac{c_3(s-N)(s-N-1)}{2\lambda}) \int_0^\infty \frac{(\lambda y)^s}{s!} e^{-\lambda y} dH_r(y) -$$

$$- c_3(s-N)\int_0^\infty \frac{\lambda^s}{s!} \int_x^\infty (z-x)^s e^{-\lambda z} dH_r(z) dx. \quad (34)$$

8.1.5 Consider the case when $s = N + N_0$.

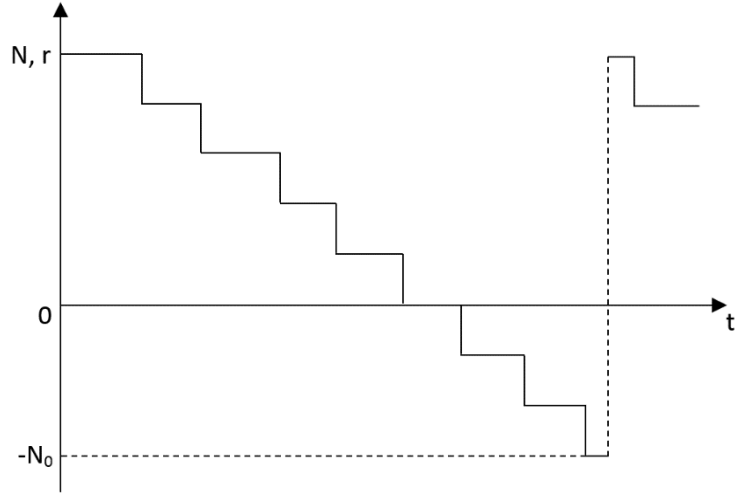

Fig. 7. Graphic illustration for the case 8.1.5

By analogy with the previously considered cases, the costs of the selling of the event $A_{N+N0}$ are composed of the product storage costs, the costs associated with acquiring a new batch of the product (replenishment) and the costs (penalties) associated with the deficit of the product. In this case, the costs associated with the deficit, it is advisable to divide into two components: losses incurred from the time of the complete exhaustion of the real stock to the moment of receipt of the last requirement, which can be taken into account, and losses incurred from the time of receipt of the last accepted requirements before the end of the delay period and replenishment.

Thus, the conditional mathematical expectation of expenses for this case is given by the relation

$$E_r(\triangle \gamma_n^{(-)} | A_{N+N_0}) = \underbrace{\sum_{i=0}^{N-1} \frac{c_1(N-i)}{\lambda}}_{\text{product storage costs}} + \underbrace{c_2(N+N_0)}_{\substack{\text{stock} \\ \text{replenishment} \\ \text{costs}}} + \underbrace{\sum_{i=1}^{N_0-1} \frac{c_3 i}{\lambda}}_{\substack{\text{costs related to the deficit before} \\ \text{the arrival of the last order}}} + \underbrace{E_r^{(3*)}(\triangle \gamma_n^{(-)} | A_{N+N_0})}_{\substack{\text{costs related to the deficit after} \\ \text{the arrival of the last order}}}.$$

From this follows the representation for the corresponding mathematical expectation of the costs by the joint distribution with the event $A_{N+N0}$

$$E_r(\triangle \gamma_n^{(-)}; A_{N+N_0}) = E_r(\triangle \gamma_n^{(-)} | A_{N+N_0}) P(A_{N+N_0}) =$$

$$= \sum_{i=0}^{N-1} \frac{c_1(N-i)}{\lambda} P(A_{N+N_0}) + c_2(N+N_0) P(A_{N+N_0}) + P(A_{N+N_0}) \sum_{i=1}^{N_0-1} \frac{c_3 i}{\lambda} + E_r^{(3*)}(\triangle \gamma_n^{(-)}; A_{N+N_0}).$$

Using Theorem 3 (formula (23))

$$E_r^{(3*)}(\triangle\gamma_n^{(-)}; A_{N+N_0}) = c_3 N_0 \tau_{r,N+N_0} = c_3 N_0 \int_0^\infty \frac{\lambda^{N+N_0}}{(N+N_0)!} \int_x^\infty (z-x)^{N+N_0} e^{-\lambda z} dH_r(z) dx.$$

According to the formula of the sum of arithmetic progression

$$\sum_{i=0}^{N-1} \frac{c_1(N-i)}{\lambda} = \frac{c_1 N(N+1)}{2\lambda}, \text{ и } \sum_{i=1}^{N_0-1} \frac{c_3 i}{\lambda} = \frac{c_3 N_0(N_0-1)}{2\lambda}.$$

Considering the above remarks, we obtain the formula for the mathematical expectation of costs by the joint distribution with the event $A_{N+N0}$

$$E_r(\triangle\gamma_n^{(-)}; A_{N+N_0}) = (\frac{c_1 N(N+1)}{2\lambda} + c_2(N+N_0) + \frac{c_3 N_0(N_0-1)}{2\lambda})\int_0^\infty \frac{(\lambda y)^{N+N_0}}{(N+N_0)!} e^{-\lambda y} dH_r(y) +$$

$$+ c_3 N_0 \int_0^\infty \frac{\lambda^{N+N_0}}{(N+N_0)!} \int_x^\infty (z-x)^{N+N_0} e^{-\lambda z} dH_r(z) dx.$$

При условии, что на данном периоде регенерации реализуется событие $A_{N+N0}$, доход образуется за счет оплаты потребителями $N + N_0$ единиц продукта, включая предварительную оплату за поставки в форме отложенного спроса. Таким образом,

Provided that the event $A_{N+N0}$ occurs on the regeneration period, the income is generated by the payment by the consumers of $N + N_0$ product units, including advance payment for deliveries in the form of deferred demand. Thus,

$$E_r(\triangle\gamma_n^{(+)} | A_{N+N_0}) = c_0(N+N_0),$$

and

$$E_r(\triangle\gamma_n^{(+)}; A_{N+N_0}) = E_r(\triangle\gamma_n^{(+)} | A_{N+N_0}) P(A_{N+N_0}) = c_0(N+N_0) \int_0^\infty \frac{(\lambda y)^{N+N_0}}{(N+N_0)!} e^{-\lambda y} dH_r(y).$$

Using the above preliminary results of the analysis of the structure of profit for the considered case, we obtain

$$E_r(\triangle\gamma_n; A_{N+N_0}) = E_r(\triangle\gamma_n^{(+)}; A_{N+N_0}) - E_r(\triangle\gamma_n^{(-)}; A_{N+N_0}) =$$
$$= (c_0(N+N_0) - \frac{c_1 N(N+1)}{2\lambda} - c_2(N+N_0) - \frac{c_3 N_0(N_0-1)}{2\lambda}) \int_0^\infty \frac{(\lambda y)^{N+N_0}}{(N+N_0)!} e^{-\lambda y} dH_r(y) -$$
$$- c_3 N_0 \int_0^\infty \frac{\lambda^{N+N_0}}{(N+N_0)!} \int_x^\infty (z-x)^{N+N_0} e^{-\lambda z} dH_r(z) dx. \tag{35}$$

8.1.6 Consider the case when $s > N + N_0$.

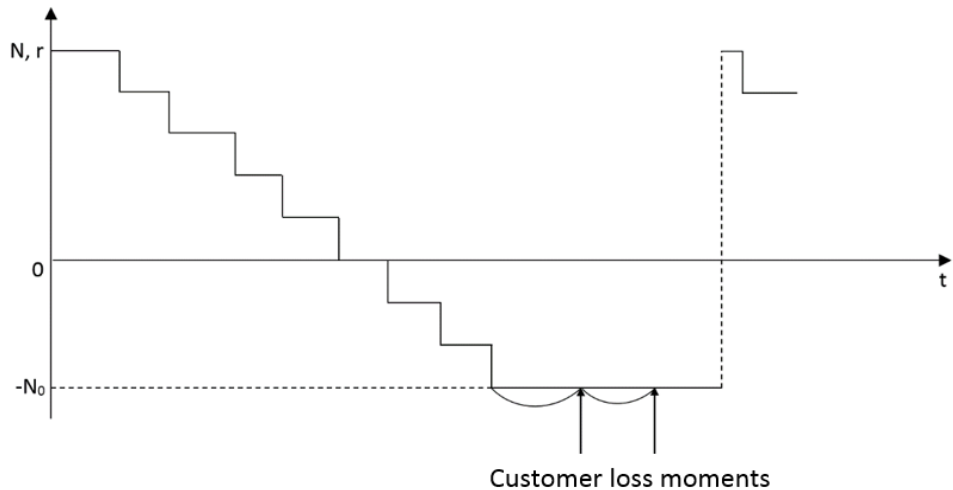

Fig. 8. Graphic illustration for the case 8.1.6

The cost structure in this case partially coincides with the structure of the corresponding value provided that the event $A_{N+N0}$ occurs. However, now part of the costs associated directly with the deficit is becoming more complicated. The indicated costs consist of the following three components:

1) the costs accumulated over a period of time from the moment the warehouse is completely empty to the moment the maximum allowable deficit of the level $N_0$ is formed;

2) costs incurred over a period of time from the moment the deficit level $N_0$ is reached until the last demand arrives (during this period the deficit remains at the level of $N_0$);

3) costs incurred over a period of time from the moment the last demand arrives until the end of the delivery delay period (in this time interval, the deficit also remains at $N_0$).

In addition, a component appears in the cost structure due to the loss of customers who cannot be taken into account. The number of such irretrievably lost customers is equal $s-(N+N_0)$. Thus, the representation for the conditional mathematical expectation of costs under the condition, that the event $A_s$, $s > N + N_0$ occurs, has the form

$$E_r(\triangle\gamma_n^{(-)} | A_s) = \underbrace{\sum_{i=0}^{N-1} \frac{c_1(N-i)}{\lambda}}_{\text{product storage costs}} + \underbrace{c_2(N+N_0)}_{\substack{\text{stock} \\ \text{replenishment} \\ \text{costs}}} + \underbrace{\sum_{i=1}^{N_0-1} \frac{c_3 i}{\lambda}}_{\substack{\text{costs related to} \\ \text{the deficit before} \\ \text{order } N+N_0}} + \underbrace{\frac{c_3 N_0}{\lambda}(s-N-N_0)}_{\substack{\text{costs related to the deficit} \\ \text{afetr order } N+N_0 \text{ but} \\ \text{before the last order}}} +$$

$$+ \underbrace{E_r^{(3*)}(\triangle\gamma_n^{(-)} | A_s)}_{\substack{\text{costs related to the deficit} \\ \text{after the arrival of the last order}}} + \underbrace{c_4^{(s-N-N_0)}}_{\substack{\text{costs of the loss} \\ \text{of customers}}}.$$

The corresponding mathematical expectation of costs, determined by the joint distribution with the event $A_s$, $s > N + N_0$, is given by the following relation, obtained taking into account the analytical transformations

$$E_r(\triangle\gamma_n^{(-)}; A_s) = E_r(\triangle\gamma_n^{(-)} | A_s)P(A_s) = (\sum_{i=0}^{N-1}\frac{c_1(N-i)}{\lambda} + c_2(N+N_0) + \sum_{i=1}^{N_0-1}\frac{c_3 i}{\lambda} +$$

$$+\frac{c_3 N_0}{\lambda}(s-N-N_0) + E_r^{(3*)}(\triangle\gamma_n^{(-)} | A_s) + c_4^{(s-N-N_0)})P(A_s) =$$

$$= \left[\frac{c_1 N(N+1)}{2\lambda} + c_2(N+N_0) + \frac{c_3 N_0(N_0-1)}{2\lambda} + \frac{c_3 N_0}{\lambda}(s-N-N_0) + c_4^{(s-N-N_0)}\right]P(A_s) +$$

$$+E_{3*}(\triangle\gamma_n^{(-)}; A_s) = (\frac{c_1 N(N+1)}{2\lambda} + c_2(N+N_0) + \frac{c_3 N_0(2s-2N-N_0-1)}{2\lambda} +$$

$$+c_4^{(s-N-N_0)})P(A_s) + E_r^{(3*)}(\triangle\gamma_n^{(-)}; A_s).$$

In the obtained expression, the value $E_r^{(3*)}(\triangle\gamma_n^{(-)}; A_s)$ is the mathematical expectation of that part of the costs that are generated in the interval from the moment the last request arrives until the end of the delay period. As in similar expressions, this mathematical expectation is determined by the joint distribution with the event $A_s$. The analytical representation of the mathematical expectation $E_r^{(3*)}(\triangle\gamma_n^{(-)}; A_s)$ is found using the statement of Theorem 3 (formula (23))

$$E_r^{(3*)}(\triangle\gamma_n^{(-)}; A_s) = c_3 N_0 \tau_{r,s} = c_3 N_0 \int_0^\infty \frac{\lambda^s}{s!}\int_x^\infty (z-x)^s e^{-\lambda z} dH_r(z) dx.$$

Given the well-known formula for the probability of receipt of demands,

$$P(A_s) = \int_0^\infty \frac{(\lambda y)^s}{s!} e^{-\lambda y} dH_r(y),$$

we obtain the final representation for the mathematical expectation of expenses by joint distribution with the event $A_s$

$$E_r(\Delta\gamma_n^{(-)}; A_s) = \left[\frac{c_1 N(N+1)}{2\lambda} + c_2(N+N_0) + \frac{c_3 N_0(2s - 2N - N_0 - 1)}{2\lambda} + \right.$$

$$\left. + c_4^{(s-N-N_0)}\right] \int_0^\infty \frac{(\lambda y)^s}{s!} e^{-\lambda y} dH_r(y) + c_3 N_0 \int_0^\infty \frac{\lambda^s}{s!} \int_x^\infty (z-x)^s e^{-\lambda z} dH_r(z) dx.$$

Since in this case the maximum allowable number of orders $N + N_0$ will be received and executed, the mathematical expectation of income by joint distribution with the event $A_s$ does not depend on the value $s > N + N_0$ and has the form

$$E_r(\Delta\gamma_n^{(+)}; A_s) = E_r(\Delta\gamma_n^{(+)} | A_s) P(A_s) = c_0(N+N_0) \int_0^\infty \frac{(\lambda y)^s}{s!} e^{-\lambda y} dH_r(y).$$

From the above results, an explicit representation follows for the mathematical expectation of profit by the joint distribution with the event $A_s$, $s > N + N_0$

$$E_r(\Delta\gamma_n; A_s) = E_r(\Delta\gamma_n^{(+)}; A_s) - E_r(\Delta\gamma_n^{(-)}; A_s) =$$

$$= c_0(N+N_0) \int_0^\infty \frac{(\lambda y)^s}{s!} e^{-\lambda y} dH_r(y) - (\frac{c_1 N(N+1)}{2\lambda} + c_2(N+N_0) + \frac{c_3 N_0(2s - 2N - N_0 - 1)}{2\lambda} +$$

$$+ c_4^{(s-N-N_0)}) \int_0^\infty \frac{(\lambda y)^s}{s!} e^{-\lambda y} dH_r(y) + c_3 N_0 \int_0^\infty \frac{\lambda^s}{s!} \int_x^\infty (z-x)^s e^{-\lambda z} dH_r(z) dx =$$

$$= (c_0(N+N_0) - \frac{c_1 N(N+1)}{2\lambda} - c_2(N+N_0) - \frac{c_3 N_0(2s - 2N - N_0 - 1)}{2\lambda} -$$

$$- c_4^{(s-N-N_0)}) \int_0^\infty \frac{(\lambda y)^s}{s!} e^{-\lambda y} dH_r(y) - c_3 N_0 \int_0^\infty \frac{\lambda^s}{s!} \int_x^\infty (z-x)^s e^{-\lambda z} dH_r(z) dx. \qquad (36)$$

We obtain the expression for the mathematical expectation $E_r(\Delta\gamma_n)$ for $r = N$. To do this, we substitute the above found representations (31) - (36) into the general relation (27). For the convenience of perceiving such a cumbersome formula, we emphasize the components in it that correspond to the considered cases for the values of the parameter s.

$$E_r(\Delta\gamma_n) = \underbrace{-c_1 N \int_0^\infty (1 - H_r(y)) dy \int_0^\infty e^{-\lambda y} dH_r(y)}_{\text{for } s=0} + \underbrace{\sum_{i=1}^{N-1} ((c_0 i - \frac{c_1 i(2N - i + 1)}{2\lambda}}_{\text{for } 1 \le s < N} -$$

$$\underbrace{- c_2 i) \int_0^\infty \frac{(\lambda y)^i}{i!} e^{-\lambda y} dH_r(y) - c_1 (N - i) \int_0^\infty \frac{\lambda^i}{i!} \int_x^\infty (z - x)^i e^{-\lambda z} dH_r(z) dx)}_{\text{for } 1 \le s < N} +$$

$$+ \underbrace{(c_0 N - \frac{c_1 N(N+1)}{2\lambda} - c_2 N) \int_0^\infty \frac{(\lambda y)^N}{N!} e^{-\lambda y} dH_r(y)}_{\text{for } s=N} + \underbrace{\sum_{i=N+1}^{N+N_0-1} ((c_0 i - \frac{c_1 N(N+1)}{2\lambda} - c_2 i -}_{\text{for } N < s < N+N_0}$$

$$\underbrace{- \frac{c_3 (i - N)(i - N - 1)}{2\lambda}) \int_0^\infty \frac{(\lambda y)^i}{i!} e^{-\lambda y} dH_r(y) - c_3 (i - N) \int_0^\infty \frac{\lambda^i}{i!} \int_x^\infty (z - x)^i e^{-\lambda z} dH_r(z) dx)}_{\text{for } N < s < N+N_0} +$$

$$+ \underbrace{(c_0 (N + N_0) - \frac{c_1 N(N+1)}{2\lambda} - c_2 (N + N_0) - \frac{c_3 N_0 (N_0 - 1)}{2\lambda}) \int_0^\infty \frac{(\lambda y)^{N+N_0}}{(N+N_0)!} e^{-\lambda y} dH_r(y) -}_{\text{for } s = N+N_0}$$

$$\underbrace{- c_3 N_0 \int_0^\infty \frac{\lambda^{N+N_0}}{(N+N_0)!} \int_x^\infty (z - x)^{N+N_0} e^{-\lambda z} dH_r(z) dx}_{\text{for } s=N+N_0} + \underbrace{\sum_{i=N+N_0+1}^\infty ((c_0 (N + N_0) - \frac{c_1 N(N+1)}{2\lambda} - c_2 (N + N_0) -}_{\text{for } s > N+N_0}$$

$$\underbrace{- \frac{c_3 N_0 (2i - 2N - N_0 - 1)}{2\lambda} - c_4^{(i-N-N_0)}) \int_0^\infty \frac{(\lambda y)^i}{i!} e^{-\lambda y} dH_r(y) - c_3 N_0 \int_0^\infty \frac{\lambda^i}{i!} \int_x^\infty (z - x)^i e^{-\lambda z} dH_r(z) dx)}_{\text{for } s > N+N_0}.$$

Now we can get the final representation for the function $C(r) = \frac{A(r)}{B(r)}$, whose global maximum determines the optimal control (see Section 6). Since, in accordance with (12) and (13),

$$A(r) = E_r(\Delta\gamma_n); \quad B(r) = E_r(\Delta t_n) = \int_0^\infty (1 - H_r(y)) dy,$$

the formula for $C(r) = I_r = \frac{E_r(\Delta\gamma_n)}{E_r(\Delta t_n)}$ has the following form

$$I_r = \frac{1}{\int_0^\infty (1-H_r(y))dy}(-c_1 N \int_0^\infty e^{-\lambda y} dH_r(y) \int_0^\infty (1-H_r(y))dy + \sum_{i=1}^{N-1}((c_0 i - \frac{c_1 i(2N-i+1)}{2\lambda} -$$

$$-c_2 i)\int_0^\infty \frac{(\lambda y)^i}{i!} e^{-\lambda y} dH_r(y) - c_1(N-i)\int_0^\infty \frac{\lambda^i}{i!}\int_x^\infty (z-x)^i e^{-\lambda z} dH_r(z)dx) +$$

$$+(c_0 N - \frac{c_1 N(N+1)}{2\lambda} - c_2 N)\int_0^\infty \frac{(\lambda y)^N}{N!} e^{-\lambda y} dH_r(y) + \sum_{i=N+1}^{N+N_0-1}((c_0 i - \frac{c_1 N(N+1)}{2\lambda} - c_2 i -$$

$$- \frac{c_3(i-N)(i-N-1)}{2\lambda})\int_0^\infty \frac{(\lambda y)^i}{i!} e^{-\lambda y} dH_r(y) - c_3(i-N)\int_0^\infty \frac{\lambda^i}{i!}\int_x^\infty (z-x)^i e^{-\lambda z} dH_r(z)dx) +$$

$$+(c_0(N+N_0) - \frac{c_1 N(N+1)}{2\lambda} - c_2(N+N_0) - \frac{c_3 N_0(N_0-1)}{2\lambda})\int_0^\infty \frac{(\lambda y)^{N+N_0}}{(N+N_0)!} e^{-\lambda y} dH_r(y) -$$

$$-c_3 N_0 \int_0^\infty \frac{\lambda^{N+N_0}}{(N+N_0)!}\int_x^\infty (z-x)^{N+N_0} e^{-\lambda z} dH_r(z)dx + \sum_{i=N+N_0+1}^\infty ((c_0(N+N_0) - \frac{c_1 N(N+1)}{2\lambda} - c_2(N+N_0) -$$

$$-\frac{c_3 N_0(2i-2N-N_0-1)}{2\lambda} - c_4^{(i-N-N_0)})\int_0^\infty \frac{(\lambda y)^i}{i!} e^{-\lambda y} dH_r(y) - c_3 N_0 \int_0^\infty \frac{\lambda^i}{i!}\int_x^\infty (z-x)^i e^{-\lambda z} dH_r(z)dx)). \quad (37)$$

We also note that formula (37) is valid for r = N. In order to preserve the unity of notation in the future, when representing the function $I_r$, analogic notation will be used.

## 8.2 The variant when $1 \leq r < N$

### 8.2.1 Consider the case when s = 0.

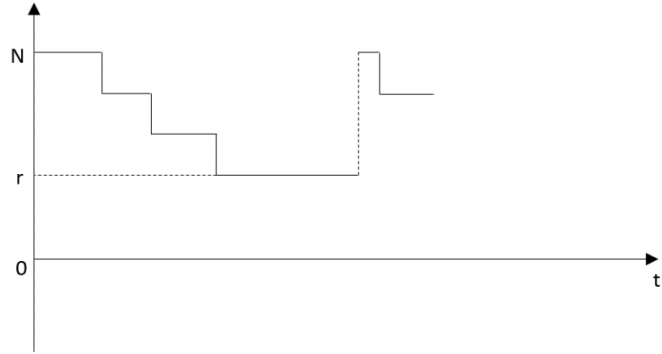

Fig. 9. Graphic illustration for the case 8.2.1

The income over the entire regeneration period in this case consists of the payment by consumers of N-r units of goods. It follows that the mathematical expectation of income by the joint distribution with event $A_0$ is determined by the formula

$$E_r(\triangle \gamma_n^{(+)}; A_0) = E_r(\triangle \gamma_n^{(+)} | A_0) P(A_0) = c_0(N-r)\int_0^\infty e^{-\lambda y} dH_r(y).$$

The costs for the considered period of regeneration consist of the product srorage costs and the stock replenishment costs. In the total amount of storage costs, it is advisable to distinguish two components: costs arising during the time interval from the beginning of the regeneration period to the beginning of the delivery delay, and costs arising during the delivery delay period. Then the conditional expectation of expenses on the regeneration period is presented in the form

$$E_r(\Delta\gamma_n^{(-)} | A_0) = \underbrace{\sum_{i=r+1}^{N} \frac{c_1 i}{\lambda}}_{\substack{product\ storage \\ costs\ before\ the\ start \\ of\ delay\ period}} + \underbrace{c_1 r \int_0^\infty (1 - H_r(y)) dy}_{\substack{product\ storage \\ costs\ after\ the\ start \\ of\ delay\ period}} + \underbrace{c_2(N-r)}_{\substack{stock \\ replenishment \\ costs}}.$$

Obviously, $\sum_{i=r+1}^{N} \frac{c_1 i}{\lambda} = \frac{c_1(N+r+1)(N-r)}{2\lambda}$, therefore, the mathematical expectation of the costs by the joint distribution with the event $A_0$ is determined by the equality

$$E_r(\Delta\gamma_n^{(-)}; A_0) = E_r(\Delta\gamma_n^{(-)} | A_0) P(A_0) = \left[ \frac{c_1(N+r+1)(N-r)}{2\lambda} + c_1 r \int_0^\infty (1 - H_r(y)) dy + \right.$$

$$\left. + c_2(N-r) \right] \int_0^\infty e^{-\lambda y} dH_r(y).$$

Now it is possible to write out explicitly the corresponding mathematical expectation of profit for the analyzed case.

$$E_r(\Delta\gamma_n; A_0) = E_r(\Delta\gamma_n^{(+)}; A_0) - E_r(\Delta\gamma_n^{(-)}; A_0) = c_0(N-r) \int_0^\infty e^{-\lambda y} dH_r(y) -$$

$$- \left[ \frac{c_1(N+r+1)(N-r)}{2\lambda} + c_1 r \int_0^\infty (1 - H_r(y)) dy + c_2(N-r) \right] \int_0^\infty e^{-\lambda y} dH_r(y) =$$

$$= \left[ c_0(N-r) - \frac{c_1(N+r+1)(N-r)}{2\lambda} - c_1 r \int_0^\infty (1 - H_r(y)) dy - c_2(N-r) \right] \int_0^\infty e^{-\lambda y} dH_r(y). \quad (38)$$

8.2.2 Consider the case when $1 \leq s < r$.

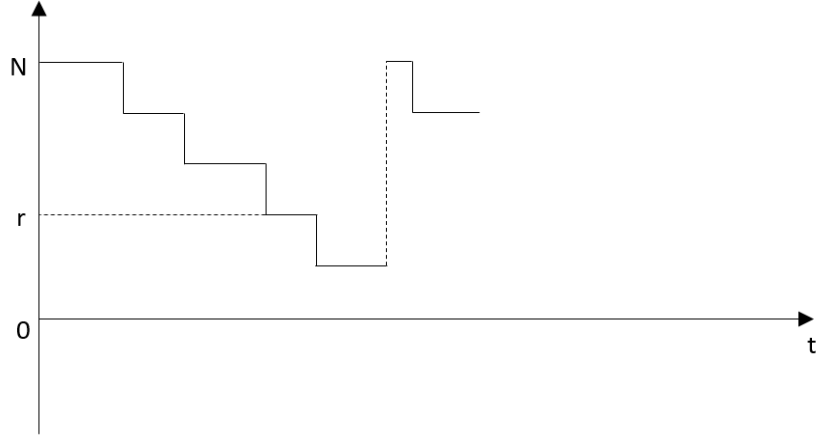

Fig. 10. Graphic illustration for the case 8.2.2

First of all, we will find the mathematical expectation of income by the joint distribution with the $A_s$ event. In this case, during the regeneration period, there will be $N - r + s$ requirements for the goods that will be met. Thus, the mathematical expectation of income by the joint distribution with the event $A_s$ has the form

$$E_r(\triangle\gamma_n^{(+)}; A_s) = E_r(\triangle\gamma_n^{(+)} | A)P(A_s) = c_0(N - r + s)\int_0^\infty \frac{(\lambda y)^s}{s!} e^{-\lambda y} dH_r(y).$$

Now we need to find the appropriate mathematical expectation of expenses during the regeneration period. As in the previous case, the costs on the regeneration period are the sum of the product storage costs and the stock replenishment costs. However, in this case, it is convenient to distinguish two other components in the total amount of storage costs: costs incurred in the time interval from the beginning of the regeneration period until the last request arrives, and costs arising from the time the last request arrives until the end of the regeneration period. Then

$$E_r(\triangle\gamma_n^{(-)}; A_s) = E_r(\triangle\gamma_n^{(-)} | A_s)P(A_s) = \left[ \underbrace{\sum_{i=r-s+1}^{N} \frac{c_1 i}{\lambda}}_{\substack{product\ storage \\ costs\ before\ the\ last \\ order}} + \underbrace{E_r^{(1*)}(\triangle\gamma_n^{(-)} | A_s)}_{\substack{product\ storage \\ costs\ after\ the\ last \\ order}} + \right.$$

$$\left. + \underbrace{c_2(N - r + s)}_{\substack{stock \\ replenishment \\ costs}} \right] P(A_s).$$

According to the formula of the sum of arithmetic progression $\sum_{i=r-s+1}^{N} \frac{c_1 i}{\lambda} = \frac{c_1(N-r+s)(N+r-s+1)}{2\lambda}$. Applying formula (23) of Theorem 3 to find the mathematical expectation $E_r^{(1*)}(\triangle\gamma_n^{(-)}; A_s) = E_r^{(1*)}(\triangle\gamma_n^{(-)} | A_s)P(A_s)$, we obtain

$$E_r^{(1*)}(\triangle\gamma_n^{(-)}; A_s) = c_1(r-s)\tau_{r,s} = c_1(r-s)\int_0^\infty \frac{\lambda^s}{s!}\int_x^\infty (z-x)^s e^{-\lambda z} dH_r(z)dx.$$

Thus, the mathematical expectation of expenses by joint distribution with the event $A_s$ is determined by the formula:

$$E_r(\triangle\gamma_n^{(-)}; A_s) = (\frac{c_1(N-r+s)(N+r-s+1)}{2\lambda} + c_2(N-r+s))\int_0^\infty \frac{(\lambda y)^s}{s!} e^{-\lambda y} dH_r(y) +$$

$$+ c_1(r-s)\int_0^\infty \frac{\lambda^s}{s!}\int_x^\infty (z-x)^s e^{-\lambda z} dH_r(z)dx.$$

Now, to obtain an explicit analytical representation for the mathematical expectation of profit by the joint distribution with the event $A_s$, it is necessary to use relation (28)

$$E_r(\triangle\gamma_n; A_s) = E_r(\triangle\gamma_n^{(+)}; A_s) - E_r(\triangle\gamma_n^{(-)}; A_s) = c_0(N-r+s)\int_0^\infty \frac{(\lambda y)^s}{s!} e^{-\lambda y} dH_r(y) -$$

$$-\left[\frac{c_1(N-r+s)(N+r-s+1)}{2\lambda} + c_2(N-r+s)\right]\int_0^\infty \frac{(\lambda y)^s}{s!} e^{-\lambda y} dH_r(y) -$$

$$-c_1(r-s)\int_0^\infty \frac{\lambda^s}{s!}\int_x^\infty (z-x)^s e^{-\lambda z} dH_r(z)dx =$$

$$= \left[c_0(N-r+s) - \frac{c_1(N-r+s)(N+r-s+1)}{2\lambda} - \right.$$

$$\left. -c_2(N-r+s)\right]\int_0^\infty \frac{(\lambda y)^s}{s!} e^{-\lambda y} dH_r(y) - c_1(r-s)\int_0^\infty \frac{\lambda^s}{s!}\int_x^\infty (z-x)^s e^{-\lambda z} dH_r(z)dx. \quad (39)$$

8.2.3 Consider the case when s = r.

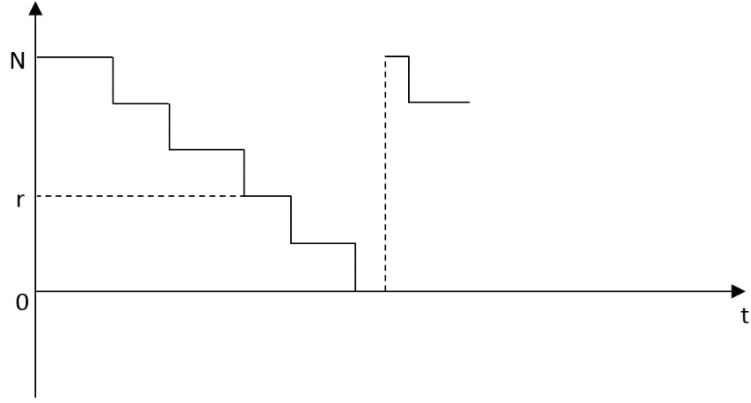

Fig. 11. Graphic illustration for the case 8.2.3

In this case, all goods that were in the warehouse at the time the regeneration period began, that is, N units, are sold. Thus, the mathematical expectation of profit by the joint distribution with the $A_r$ event is determined by the formula:

$$E_r(\triangle\gamma_n^{(+)}; A_r) = E_r(\triangle\gamma_n^{(+)} | A_r)P(A_r) = c_0 N \int_0^\infty \frac{(\lambda y)^r}{r!} e^{-\lambda y} dH_r(y).$$

As before, the costs during the regeneration period consist of the product storage costs and the stock replenishment costs. However, provided that the event $A_r$ occurs, the entire stock of goods will be used up until the end of the regeneration period. Then the costs associated with storing unspent product are equal to zero.

$$E_r(\triangle\gamma_n^{(-)} | A_r) = \underbrace{\sum_{i=1}^{N} \frac{c_1 i}{\lambda}}_{\substack{\text{product} \\ \text{storage} \\ \text{costs}}} + \underbrace{c_2 N}_{\substack{\text{stock} \\ \text{replenishment} \\ \text{costs}}} = \frac{c_1(N+1)N}{2\lambda} + c_2 N.$$

Thus,

$$E_r(\triangle\gamma_n^{(-)}; A_r) = E_r(\triangle\gamma_n^{(-)} | A_r)P(A_r) = \left[\frac{c_1(N+1)N}{2\lambda} + c_2 N\right]P(A_r) =$$

$$= \left[\frac{c_1(N+1)N}{2\lambda} + c_2 N\right] \int_0^\infty \frac{(\lambda y)^r}{r!} e^{-\lambda y} dH_r(y).$$

Therefore, the mathematical expectation of profit by the joint distribution with the event $A_r$ has the form

$$E_r(\triangle\gamma_n; A_r) = E_r(\triangle\gamma_n^{(+)}; A_r) - E_r(\triangle\gamma_n^{(-)}; A_r) = c_0 N \int_0^\infty \frac{(\lambda y)^r}{r!} e^{-\lambda y} dH_r(y) -$$

$$-(\frac{c_1(N+1)N}{2\lambda} + c_2 N)\int_0^\infty \frac{(\lambda y)^r}{r!} e^{-\lambda y} dH_r(y) =$$

$$= (c_0 N - \frac{c_1(N+1)N}{2\lambda} - c_2 N)\int_0^\infty \frac{(\lambda y)^r}{r!} e^{-\lambda y} dH_r(y). \qquad (40)$$

8.2.4 Consider the case when $r < s < r + N_0$.

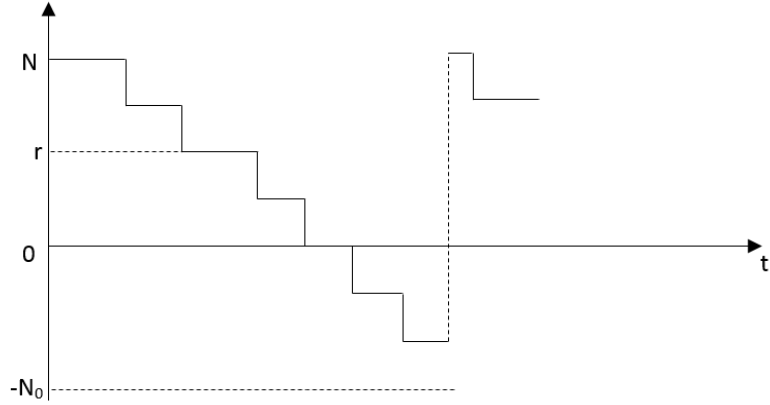

Fig. 12. Graphic illustration for the case 8.2.4

First, we find an explicit representation for the mathematical expectation of income by the joint distribution with the event $A_s$. Note that in this case $N - r + s$ units of goods will be sold, including $N-r$ units before the start of the delivery delay period and $s$ units during the delay period. In the standard way we obtain

$$E_r(\triangle\gamma_n^{(+)}; A_s) = E_r(\triangle\gamma_n^{(+)} | A_s) P(A_s) = c_0(N - r + s)\int_0^\infty \frac{(\lambda y)^s}{s!} e^{-\lambda y} dH_r(y).$$

The costs consist of the following components: the product storage costs, the stock replenishment costs and the costs (penalties) associated with the deficit. In this case, penalties associated with the deficit should be divided into penalties that arise over the period from the moment the deficit is formed until the last request arrives, and penalties that arise over the period from the moment the last request arrives until the end of the regeneration period. We note in addition that the moment of deficit formation coincides with the moment of receipt of the demand with the number $N + 1$, counting from the beginning of the regeneration period, at this moment the deficit takes the value 1. The last requirement for the regeneration period has the number $N - r + s$, at the moment the deficit reaches the level of s - r and then does not change until the end of the regeneration period. From here follows the representation for the conditional mathematical expectation of expenses under the condition that the event $A_s$ occurs:

$$E_r(\triangle \gamma_n^{(-)} | A_s) = \sum_{i=1}^{N} \frac{c_1 i}{\lambda} + \underbrace{c_2(N-r+s)}_{\substack{\text{stock} \\ \text{replenishment} \\ \text{costs}}} + \underbrace{\sum_{i=1}^{s-r-1} \frac{c_3 i}{\lambda}}_{\substack{\text{costs related to the deficit} \\ \text{before the arrival of} \\ \text{the last order}}} + \underbrace{E_r^{(3*)}(\triangle \gamma_n^{(-)} | A_s)}_{\substack{\text{costs related to the deficit} \\ \text{after the arrival of} \\ \text{the last order}}}.$$

product storage costs

Using the formula of the sum of arithmetic progression

$$\sum_{i=1}^{N} \frac{c_1 i}{\lambda} = \frac{c_1(N+1)N}{2\lambda}, \quad \sum_{i=1}^{s-r-1} \frac{c_3 i}{\lambda} = \frac{c_3(s-r)(s-r-1)}{2\lambda}.$$

Passing from the conditional mathematical expectation to the mathematical expectation by the joint distribution, we obtain:

$$E_r(\triangle \gamma_n^{(-)}; A_s) = E_r(\triangle \gamma_n^{(-)} | A_s) P(A_s) = \left[ \frac{c_1(N+1)N}{2\lambda} + c_2(N-r+s) + \right.$$

$$\left. + \frac{c_3(s-r)(s-r-1)}{2\lambda} \right] P(A_s) + E_r^{(3*)}(\triangle \gamma_n^{(-)}; A_s).$$

We express the value $E_r^{(3*)}(\triangle \gamma_n^{(-)}; A_s)$ using Theorem 3

$$E_r^{(3*)}(\triangle \gamma_n^{(-)}; A_s) = c_3(s-r)\tau_{r,s} = c_3(s-r)\int_0^\infty \frac{\lambda^s}{s!} \int_x^\infty (z-x)^s e^{-\lambda z} dH_r(z) dx.$$

Then,

$$E_r(\triangle \gamma_n^{(-)}; A_s) = (\frac{c_1(N+1)N}{2\lambda} + c_2(N-r+s) + + \frac{c_3(s-r)(s-r-1)}{2\lambda}) \int_0^\infty \frac{(\lambda y)^s}{s!} e^{-\lambda y} dH_r(y) +$$

$$+ c_3(s-r) \int_0^\infty \frac{\lambda^s}{s!} \int_x^\infty (z-x)^s e^{-\lambda z} dH_r(z) dx.$$

From the results obtained there follows a representation for the mathematical expectation of profit on a certain regeneration period, determined by the joint distribution with the event $A_s$:

$$E_r(\Delta\gamma_n; A_s) = E_r(\Delta\gamma_n^{(+)}; A_s) - E_r(\Delta\gamma_n^{(-)}; A_s) = c_0(N-r+s)\int_0^\infty \frac{(\lambda y)^s}{s!} e^{-\lambda y} dH_r(y) -$$

$$-(\frac{c_1(N+1)N}{2\lambda} + c_2(N-r+s) + \frac{c_3(s-r)(s-r-1)}{2\lambda})\int_0^\infty \frac{(\lambda y)^s}{s!} e^{-\lambda y} dH_r(y) -$$

$$-c_3(s-r)\int_0^\infty \frac{\lambda^s}{s!}\int_x^\infty (z-x)^s e^{-\lambda z} dH_r(z)dx = (c_0(N-r+s) - \frac{c_1(N+1)N}{2\lambda} - c_2(N-r+s)$$

$$- \frac{c_3(s-r)(s-r-1)}{2\lambda})\int_0^\infty \frac{(\lambda y)^s}{s!} e^{-\lambda y} dH_r(y) - c_3(s-r)\int_0^\infty \frac{\lambda^s}{s!}\int_x^\infty (z-x)^s e^{-\lambda z} dH_r(z)dx. \quad (41)$$

8.2.5 Consider the case when $s = r + N_0$.

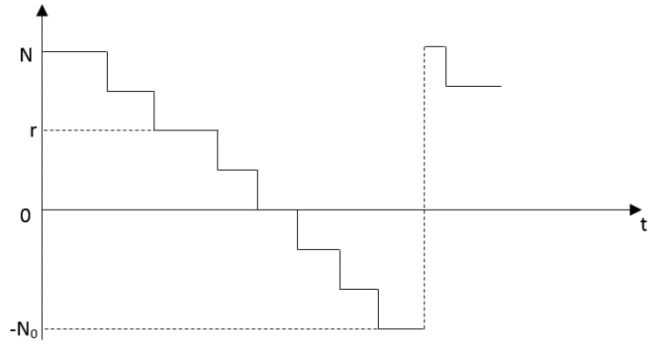

Fig. 13. Graphic illustration for the case 8.2.5

In the considered case, during the regeneration period exactly $N + N_0$ requests are received for the goods, including $N - r$ before the delay and $s = r + N_0$ during the delivery delay. In accordance with our model, all these requirements will be satisfied ($N$ units are real stock and $N_0$ units are deferred demand).

$$E_r(\Delta\gamma_n^{(+)}; A_{r+N_0}) = E_r(\Delta\gamma_n^{(+)} | A_{r+N_0}) P(A_{r+N_0}) = c_0(N+N_0)\int_0^\infty \frac{(\lambda y)^{r+N_0}}{(r+N_0)!} e^{-\lambda y} dH_r(y).$$

The costs in this case are the sum of the product storage costs, the costs necessary to replenish the stock, and the costs arising from a deficit of goods. At the same time, the costs associated with the deficit should be divided into two components: costs that arise during the period from the moment the deficit is formed until the last request arrives, and expenses that arise during the period from the moment the last request is received until the end of the regeneration period. From here follows the representation for the conditional mathematical expectation of costs, provided that the event $A_{r+N_0}$ occurs:

$$E_r(\Delta\gamma_n^{(-)} | A_{r+N_0}) = \sum_{i=1}^{N}\frac{c_1 i}{\lambda} + \underbrace{c_2(N+N_0)}_{\substack{\text{stock} \\ \text{replenishment} \\ \text{costs}}} + \underbrace{\sum_{i=1}^{N_0-1}\frac{c_3 i}{\lambda}}_{\substack{\text{costs related to the deficit} \\ \text{before the arrival of} \\ \text{the last order}}} + \underbrace{E_r^{(3^*)}(\Delta\gamma_n^{(-)} | A_{r+N_0})}_{\substack{\text{costs related to the deficit} \\ \text{after the arrival of} \\ \text{the last order}}}.$$

<span style="text-indent:0">product storage costs</span>

According to the formula of the sum of arithmetic progression

$$\sum_{i=1}^{N}\frac{c_1 i}{\lambda} = \frac{c_1(N+1)N}{2\lambda}, \text{ и } \sum_{i=1}^{N_0-1}\frac{c_3 i}{\lambda} = \frac{c_3 N_0(N_0-1)}{2\lambda}.$$

Now we turn to the mathematical expectation of the costs by joint distribution.

$$E_r(\Delta\gamma_n^{(-)}; A_{r+N_0}) = E_r(\Delta\gamma_n^{(-)} | A_{r+N_0})P(A_{r+N_0}) = \left[\frac{c_1(N+1)N}{2\lambda} + c_2(N+N_0) + \right.$$
$$\left. + \frac{c_3 N_0(N_0-1)}{2\lambda}\right]P(A_{r+N_0}) + E_r^{(3^*)}(\Delta\gamma_n^{(-)}; A_{r+N_0}).$$

Using Theorem 3 (formula (23)), we obtain

$$E_r^{(3^*)}(\Delta\gamma_n^{(-)}; A_{r+N_0}) = c_3 N_0 \tau_{r,r+N_0} = c_3 N_0 \int_0^\infty \frac{\lambda^{r+N_0}}{(r+N_0)!}\int_x^\infty (z-x)^{r+N_0} e^{-\lambda z} dH_r(z)dx.$$

Thus,

$$E_r(\Delta\gamma_n^{(-)}; A_{r+N_0}) = \left[\frac{c_1(N+1)N}{2\lambda} + c_2(N+N_0) + + \frac{c_3 N_0(N_0-1)}{2\lambda}\right]\int_0^\infty \frac{(\lambda y)^{r+N_0}}{(r+N_0)!}e^{-\lambda y}dH_r(y) +$$
$$+ c_3 N_0 \int_0^\infty \frac{\lambda^{r+N_0}}{(r+N_0)!}\int_x^\infty (z-x)^{r+N_0} e^{-\lambda z} dH_r(z)dx.$$

Now it becomes possible to write out an explicit analytical representation for the mathematical expectation of profit by the joint distribution with the event $A_{r+N0}$ for the case under consideration.

$$E_r(\Delta\gamma_n; A_{r+N_0}) = E_r(\Delta\gamma_n^{(+)}; A_{r+N_0}) - E_r(\Delta\gamma_n^{(-)}; A_{r+N_0}) = c_0(N+N_0)\int_0^\infty \frac{(\lambda y)^{r+N_0}}{r+N_0!}e^{-\lambda y}dH_r(y) -$$
$$-\left[\frac{c_1(N+1)N}{2\lambda} + c_2(N+N_0) + \frac{c_3 N_0(N_0-1)}{2\lambda}\right]\int_0^\infty \frac{(\lambda y)^{r+N_0}}{(r+N_0)!}e^{-\lambda y}dH_r(y) -$$
$$-c_3 N_0 \int_0^\infty \frac{\lambda^{r+N_0}}{r+N_0!}\int_x^\infty (z-x)^{r+N_0} e^{-\lambda z} dH_r(z)dx = \left[c_0(N+N_0) - \frac{c_1(N+1)N}{2\lambda} - c_2(N+N_0) -\right.$$
$$\left. -\frac{c_3 N_0(N_0-1)}{2\lambda}\right]\int_0^\infty \frac{(\lambda y)^{r+N_0}}{(r+N_0)!}e^{-\lambda y}dH_r(y) - c_3 N_0 \int_0^\infty \frac{\lambda^{r+N_0}}{r+N_0!}\int_x^\infty (z-x)^{r+N_0} e^{-\lambda z} dH_r(z)dx. \quad (42)$$

8.2.6 Consider the case when $s > r + N_0$.

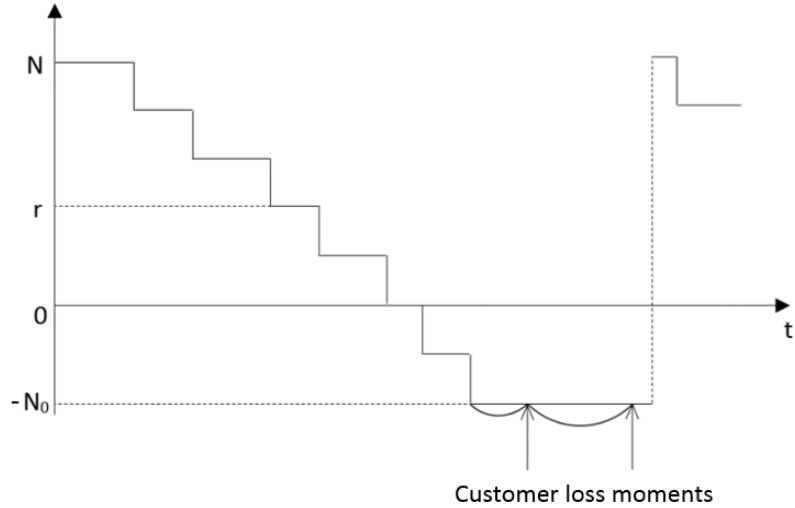

Fig. 14. Graphic illustration for the case 8.2.6

Let us first note that for any $s > r + N_0$ the maximum possible number of the product units is sold over the regeneration period, namely $N + N_0$.

$$E_r(\Delta\gamma_n^{(+)}; A_s) = E_r(\Delta\gamma_n^{(+)} | A_s)P(A_s) = c_0(N + N_0)\int_0^\infty \frac{(\lambda y)^s}{s!}e^{-\lambda y}dH_r(y).$$

In the considered variant, the costs have the most complicated structure and consist of the following components:
- real stock storage costs;
- stock replenishment costs;
- costs related to product deficit;
- costs resulting from loss of clients (orders).

At the same time, the costs related to the product deficit should be divided into three parts with different analytical expressions. The first part of these costs occur over the interval from the deficit formation point till the arrival of the last order that can be completed, i.e. the order numbered $r + N_0$ starting from the beginning of the delay period. The second part of these costs occurs over the time interval from the arrival of the order numbered $r + N_0$ till the arrival of the last order numbered $s$, also starting from the beginning of the delay period. Note that over the whole interval, the deficit has a constant value $N_0$, and the mathematical expectation of this interval equals $\frac{s-r-N_0}{\lambda}$. The third part of these costs occurs over the interval from the arrival of the last order till the end of the regeneration period. Denote the conditional mathematical expectation of the last cost deficit component by $E_r^{(3*)}(\Delta\gamma_n^{(-)} | A_s)$, similarly to the corresponding storage costs. Finally, in the variant considered, we will lose $s-r-N_0$ orders that will neither be completed really nor in the form of deferred demand. Thus, the conditional mathematical expectation of costs if event $A_s$ occurs looks as

$$E_r(\Delta\gamma_n^{(-)} | A_s) = \underbrace{\sum_{i=1}^{N} \frac{c_1 i}{\lambda}}_{\substack{product \\ storage \\ costs}} + \underbrace{c_2(N+N_0)}_{\substack{stock \\ replenishment \\ costs}} + \underbrace{\sum_{i=1}^{N_0-1} \frac{c_3 i}{\lambda}}_{\substack{costs\ related\ to \\ the\ deficit\ before \\ order\ r+N_0}} + \underbrace{\frac{c_3 N_0}{\lambda}(s-r-N_0)}_{\substack{costs\ related\ to\ the\ deficit \\ after\ order\ r+N_0\ but \\ before\ the\ last\ order}} +$$

$$+ \underbrace{E_r^{(3*)}(\Delta\gamma_n^{(-)} | A_s)}_{\substack{costs\ related\ to\ the\ deficit \\ after\ the\ arrival \\ of\ the\ last\ order}} + \underbrace{c_4^{(s-r-N_0)}}_{\substack{costs\ of\ the\ loss \\ of\ clients}}.$$

Taking into account that $\sum_{i=1}^{N} \frac{c_1 i}{\lambda} = \frac{c_1(N+1)N}{2\lambda}$ and $\sum_{i=1}^{N_0-1} \frac{c_3 i}{\lambda} = \frac{c_3 N_0(N_0-1)}{2\lambda}$, we get the following expression for the mathematical expectation of the costs by the joint distribution with event $A_s$

$$E_r(\Delta\gamma_n^{(-)}; A_s) = E_r(\Delta\gamma_n^{(-)} | A_s)P(A_s) = \left[\frac{c_1(N+1)N}{2\lambda} + c_2(N+N_0) + \right.$$

$$\left. + \frac{c_3 N_0(N_0-1)}{2\lambda} + \frac{c_3 N_0}{\lambda}(s-r-N_0) + c_4^{(s-r-N_0)}\right]P(A_s) + E_r^{(3*)}(\Delta\gamma_n^{(-)}; A_s).$$

According to the statement of theorem 3,

$$E_r^{(3*)}(\Delta\gamma_n^{(-)}; A_s) = c_3 N_0 \tau_{r,s} = c_3 N_0 \int_0^\infty \frac{\lambda^s}{s!} \int_x^\infty (z-x)^s e^{-\lambda z} dH_r(z) dx.$$

In accordance with the latter note, the formula for $E_r(\Delta\gamma_n^{(-)}; A_s)$ looks as follows:

$$E_r(\Delta\gamma_n^{(-)}; A_s) = \left[\frac{c_1(N+1)N}{2\lambda} + c_2(N+N_0) + \frac{c_3 N_0(2s-2r-N_0-1)}{2\lambda} + \right.$$

$$\left. + c_4^{(s-r-N_0)}\right]\int_0^\infty \frac{(\lambda y)^s}{s!} e^{-\lambda y} dH_r(y) + c_3 N_0 \int_0^\infty \frac{\lambda^s}{s!} \int_x^\infty (z-x)^s e^{-\lambda z} dH_r(z) dx.$$

The obtained results allow us to express the mathematical expectation of the profit by the joint distribution with event $A_s$:

$$E_r(\triangle\gamma_n; A_s) = E_r(\triangle\gamma_n^{(+)}; A_s) - E_r(\triangle\gamma_n^{(-)}; A_s) = c_0(N + N_0)\int_0^\infty \frac{(\lambda y)^s}{s!}e^{-\lambda y}dH_r(y) -$$

$$-\left[\frac{c_1(N+1)N}{2\lambda} + c_2(N + N_0) + \frac{c_3 N_0(2s - 2r - N_0 - 1)}{2\lambda} + \right.$$

$$\left. + c_4^{(s-r-N_0)}\right]\int_0^\infty \frac{(\lambda y)^s}{s!}e^{-\lambda y}dH_r(y) - c_3 N_0 \int_0^\infty \frac{\lambda^s}{s!}\int_x^\infty (z-x)^s e^{-\lambda z}dH_r(z)dx =$$

$$=\left[c_0(N + N_0) - \frac{c_1(N+1)N}{2\lambda} - c_2(N + N_0) - \frac{c_3 N_0(2s - 2r - N_0 - 1)}{2\lambda} - \right.$$

$$\left. - c_4^{(s-r-N_0)}\right]\int_0^\infty \frac{(\lambda y)^s}{s!}e^{-\lambda y}dH_r(y) - c_3 N_0 \int_0^\infty \frac{\lambda^s}{s!}\int_x^\infty (z-x)^s e^{-\lambda z}dH_r(z)dx. \quad (43)$$

Applying the obtained relations, we can complete the study of the variant in which the control parameter r satisfies the condition: $1 \leq r < N$. Let us first find an expression for the mathematical expectation of the profit over the regeneration period if the solution r is accepted. For convenience, let us emphasize the members of the obtained expression corresponding to the various ranges of values of parameter s.

$$E_r(\triangle\gamma_n) = \underbrace{\left[c_0(N-r) - \frac{c_1(N+r+1)(N-r)}{2\lambda} - c_1 r\int_0^\infty (1-H_r(y))dy - c_2(N-r))\int_0^\infty e^{-\lambda y}dH_r(y) + \right.}_{for\ s=0}$$

$$+ \underbrace{\sum_{i=1}^{r-1}((c_0(N-r+i) - \frac{c_1(N-r+i)(N+r-i+1)}{2\lambda} - c_2(N-r+i))\int_0^\infty \frac{(\lambda y)^i}{i!}e^{-\lambda y}dH_r(y) -}_{for\ 1\leq s<r}$$

$$\underbrace{-c_1(r-i)\int_0^\infty \frac{\lambda^i}{i!}\int_x^\infty (z-x)^i e^{-\lambda z}dH_r(z)dx) + (c_0 N - \frac{c_1(N+1)N}{2\lambda} - c_2 N)\int_0^\infty \frac{(\lambda y)^r}{r!}e^{-\lambda y}dH_r(y) +}_{for\ s=r}$$

$$+ \underbrace{\sum_{i=r+1}^{r+N_0-1}((c_0(N-r+i) - \frac{c_1(N+1)N}{2\lambda} - c_2(N-r+i) - \frac{c_3(i-r)(i-r-1)}{2\lambda})\int_0^\infty \frac{(\lambda y)^i}{i!}e^{-\lambda y}dH_r(y) -}_{for\ r<s<r+N_0}$$

$$\underbrace{-c_3(i-r)\int_0^\infty \frac{\lambda^i}{i!}\int_x^\infty (z-x)^i e^{-\lambda z}dH_r(z)dx) + (c_0(N+N_0) - \frac{c_1(N+1)N}{2\lambda} - c_2(N+N_0) -}_{for\ s=r+N_0}$$

$$- \underbrace{\frac{c_3 N_0(N_0-1)}{2\lambda})\int_0^\infty \frac{(\lambda y)^{r+N_0}}{(r+N_0)!}e^{-\lambda y}dH_r(y) - c_3 N_0 \int_0^\infty \frac{\lambda^{r+N_0}}{r+N_0!}\int_x^\infty (z-x)^{r+N_0}e^{-\lambda z}dH_r(z)dx +}$$

$$+ \underbrace{\sum_{i=r+N_0+1}^\infty ((c_0(N+N_0) - \frac{c_1(N+1)N}{2\lambda} - c_2(N+N_0) - \frac{c_3 N_0(2i-2r-N_0-1)}{2\lambda} -}_{for\ s>r+N_0}$$

$$\underbrace{-c_4^{(i-r-N_0)})\int_0^\infty \frac{(\lambda y)^i}{i!}e^{-\lambda y}dH_r(y) - c_3 N_0 \int_0^\infty \frac{\lambda^i}{i!}\int_x^\infty (z-x)^i e^{-\lambda z}dH_r(z)dx\right].}_{for\ s>r+N_0}$$

Now we can find the final analytical representation for the management efficiency indicator $I_r$.

$$I_r = \frac{1}{\frac{N-r}{\lambda} + \int_0^\infty (1-H_r(y))dy} \Bigg[ (c_0(N-r) - \frac{c_1(N+r+1)(N-r)}{2\lambda} - c_1 r \int_0^\infty (1-H_r(y))dy -$$

$$-c_2(N-r))\int_0^\infty e^{-\lambda y}dH_r(y) + \sum_{i=1}^{r-1}((c_0(N-r+i) - \frac{c_1(N-r+i)(N+r-i+1)}{2\lambda} - c_2(N-r+i)) *$$

$$* \int_0^\infty \frac{(\lambda y)^i}{i!} e^{-\lambda y}dH_r(y) - c_1(r-i)\int_0^\infty \frac{\lambda^i}{i!} \int_x^\infty (z-x)^i e^{-\lambda z}dH_r(z)dx) +$$

$$+ (c_0 N - \frac{c_1(N+1)N}{2\lambda} - c_2 N)\int_0^\infty \frac{(\lambda y)^r}{r!} e^{-\lambda y}dH_r(y) +$$

$$+ \sum_{i=r+1}^{r+N_0-1}((c_0(N-r+i) - \frac{c_1(N+1)N}{2\lambda} - c_2(N-r+i) - \frac{c_3(i-r)(i-r-1)}{2\lambda})\int_0^\infty \frac{(\lambda y)^i}{i!} e^{-\lambda y}dH_r(y) -$$

$$- c_3(i-r)\int_0^\infty \frac{\lambda^i}{i!} \int_x^\infty (z-x)^i e^{-\lambda z}dH_r(z)dx) + (c_0(N+N_0) - \frac{c_1(N+1)N}{2\lambda} - c_2(N+N_0) -$$

$$- \frac{c_3 N_0(N_0-1)}{2\lambda})\int_0^\infty \frac{(\lambda y)^{r+N_0}}{(r+N_0)!} e^{-\lambda y}dH_r(y) - c_3 N_0 \int_0^\infty \frac{\lambda^{r+N_0}}{r+N_0!} \int_x^\infty (z-x)^{r+N_0} e^{-\lambda z}dH_r(z)dx +$$

$$+ \sum_{i=r+N_0+1}^\infty ((c_0(N+N_0) - \frac{c_1(N+1)N}{2\lambda} - c_2(N+N_0) - \frac{c_3 N_0(2i-2r-N_0-1)}{2\lambda} -$$

$$- c_4^{(i-r-N_0)})\int_0^\infty \frac{(\lambda y)^i}{i!} e^{-\lambda y}dH_r(y) - c_3 N_0 \int_0^\infty \frac{\lambda^i}{i!} \int_x^\infty (z-x)^i e^{-\lambda z}dH_r(z)dx \Bigg]. \qquad (44)$$

## 8.3 The variant when r = 0

### 8.3.1 Consider the case when s = 0.

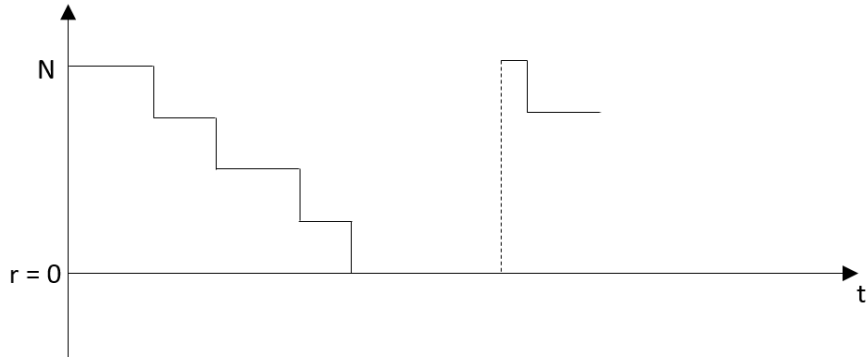

Fig. 15. Graphic illustration for the case 8.3.1

We begin the analysis of this option with the case in which during the period of delivery delays not a single client arrives. Note that if the condition r = 0 is fulfilled and the event $A_0$ occurs, then exactly N requirements will arrive during the regeneration period, all before the delivery delay beginning. Received requirements will be satisfied due to the available real stock. Also $E_r(\triangle\gamma_n^{(+)} | A_0) = c_0 N$. It follows that

$$E_r(\triangle\gamma_n^{(+)}; A_0) = E_r(\triangle\gamma_n^{(+)} | A_0)P(A_0) = c_0 N \int_0^\infty e^{-\lambda y} dH_r(y).$$

The costs in this option are the sum of the product storage costs and the costs associated with replenishment

$$E_r(\triangle\gamma_n^{(-)}; A_0) = E_r(\triangle\gamma_n^{(-)} | A_0)P(A_0) = [\underbrace{\sum_{i=1}^N \frac{c_1 i}{\lambda}}_{\substack{\text{product}\\\text{storage}\\\text{costs}}} + \underbrace{c_2 N}_{\substack{\text{stock}\\\text{replenishment}\\\text{costs}}}]\int_0^\infty e^{-\lambda y} dH_r(y) =$$

$$= \left[\frac{c_1 N(N+1)}{2\lambda} + c_2 N\right]\int_0^\infty e^{-\lambda y} dH_r(y).$$

Thus, the mathematical expectation of profit during the regeneration period by the joint distribution with event $A_0$ is determined by the following relation

$$E_r(\triangle\gamma_n; A_0) = E_r(\triangle\gamma_n^{(+)}; A_0) - E_r(\triangle\gamma_n^{(-)}; A_0) = c_0 N \int_0^\infty e^{-\lambda y} dH_r(y) -$$

$$-\left[\frac{c_1 N(N+1)}{2\lambda} + c_2 N\right]\int_0^\infty e^{-\lambda y} dH_r(y) = \left[c_0 N - \frac{c_1 N(N+1)}{2\lambda} - c_2 N\right]\int_0^\infty e^{-\lambda y} dH_r(y). \quad (45)$$

8.3.2 Consider the case when $1 \leq s < N_0$.

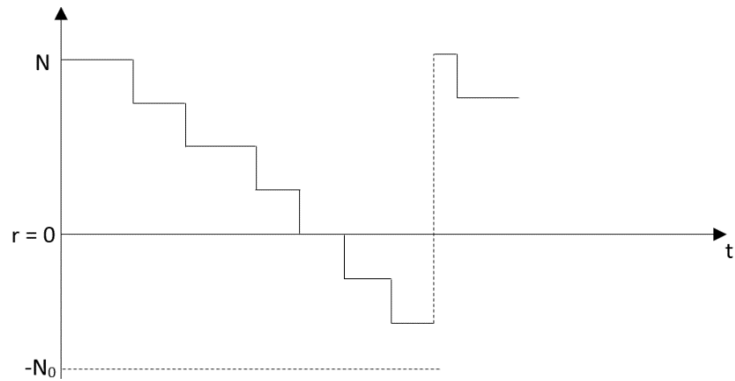

Fig. 16. Graphic illustration for the case 8.3.2

In this case, event $A_s$ occurs, so $N + s$ units of goods are sold in the system. Then

$$E_r(\Delta\gamma_n^{(+)}; A_s) = E_r(\Delta\gamma_n^{(+)} | A_s)P(A_s) = c_0(N+s)\int_0^\infty \frac{(\lambda y)^s}{s!}e^{-\lambda y}dH_r(y).$$

The costs in this option are the sum of the product storage costs, the costs necessary to replenish the stock, and the costs (penalties) associated with the deficit of goods. To replenish inventory, N + s units of goods are required. Moreover, the penalties associated with the deficit are divided into two components: penalties that arise from the moment the deficit is formed until the last request arrives, and penalties that arise from the time the last request arrives until the end of the regeneration period. By analogy with the previous results, we have

$$E_r(\Delta\gamma_n^{(-)} | A_s) = \underbrace{\sum_{i=1}^{N}\frac{c_1 i}{\lambda}}_{\substack{\text{product}\\\text{storage}\\\text{costs}}} + \underbrace{c_2(N+s)}_{\substack{\text{stock}\\\text{replenishment}\\\text{costs}}} + \underbrace{\sum_{i=1}^{s-1}\frac{c_3 i}{\lambda}}_{\substack{\text{costs related to the deficit}\\\text{before the arrival of}\\\text{the last order}}} + \underbrace{E_r^{(3*)}(\Delta\gamma_n^{(-)} | A_s)}_{\substack{\text{costs related to the deficit}\\\text{after the arrival of}\\\text{the last order}}} =$$

$$= \frac{c_1 N(N+1)}{2\lambda} + c_2(N+s) + \frac{c_3 s(s-1)}{2\lambda} + E_r^{(3*)}(\Delta\gamma_n^{(-)} | A_s).$$

Passing to the mathematical expectation of expenses by joint distribution with the event $A_s$, we obtain

$$E_r(\Delta\gamma_n^{(-)}; A_s) = E_r(\Delta\gamma_n^{(-)} | A_s)P(A_s) = (\frac{c_1 N(N+1)}{2\lambda} + c_2(N+s) + \frac{c_3 s(s-1)}{2\lambda})P(A_s) +$$
$$+ E_r^{(3*)}(\Delta\gamma_n^{(-)}; A_s).$$

Based on Theorem 3

$$E_r^{(3*)}(\Delta\gamma_n^{(-)}; A_s) = c_3 s\tau_{r,s} = c_3 s\int_0^\infty \frac{\lambda^s}{s!}\int_x^\infty (z-x)^s e^{-\lambda z}dH_r(z)dx.$$

Thus, the mathematical expectation of expenses by joint distribution with the event $A_s$ is analytically expressed by the formula

$$E_r(\Delta\gamma_n^{(-)}; A_s) = (\frac{c_1 N(N+1)}{2\lambda} + c_2(N+s) +$$
$$+ \frac{c_3 s(s-1)}{2\lambda})\int_0^\infty \frac{(\lambda y)^s}{s!}e^{-\lambda y}dH_r(y) + c_3 s\int_0^\infty \frac{\lambda^s}{s!}\int_x^\infty (z-x)^s e^{-\lambda z}dH_r(z)dx.$$

Based on the previous remarks, we get an idea for the mathematical expectation of profit on the regeneration period, determined by the joint distribution with the event $A_s$:

$$E_r(\triangle\gamma_n; A_s) = E_r(\triangle\gamma_n^{(+)}; A_s) - E_r(\triangle\gamma_n^{(-)}; A_s) = c_0(N+s)\int_0^\infty \frac{(\lambda y)^s}{s!} e^{-\lambda y} dH_r(y) -$$

$$-\left[\frac{c_1 N(N+1)}{2\lambda} + c_2(N+s) + \frac{c_3 s(s-1)}{2\lambda}\right]\int_0^\infty \frac{(\lambda y)^s}{s!} e^{-\lambda y} dH_r(y) - c_3 s\int_0^\infty \frac{\lambda^s}{s!}\int_x^\infty (z-x)^s e^{-\lambda z} dH_r(z) dx =$$

$$= \left[c_0(N+s) - \frac{c_1 N(N+1)}{2\lambda} - c_2(N+s) - \frac{c_3 s(s-1)}{2\lambda}\right]\int_0^\infty \frac{(\lambda y)^s}{s!} e^{-\lambda y} dH_r(y) -$$

$$- c_3 s\int_0^\infty \frac{\lambda^s}{s!}\int_x^\infty (z-x)^s e^{-\lambda z} dH_r(z) dx \qquad (46)$$

8.3.3 Consider the case when $s = N_0$.

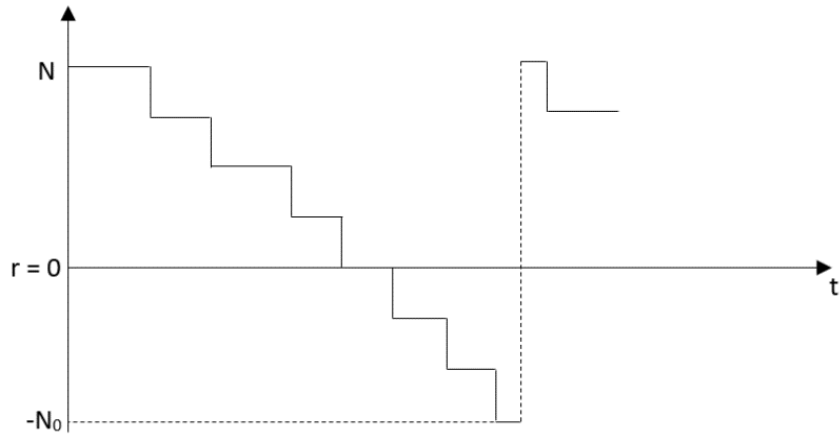

Fig. 17. Graphic illustration for the case 8.3.3

In this case, exactly $N + N_0$ customers will come into the system during the regeneration period, all of which will be satisfied (N due to real stock, $N_0$ in the form of deferred demand). Thus, the mathematical expectation of income during the regeneration period by the joint distribution with the event $A_{N0}$ is determined by the formula

$$E_r(\triangle\gamma_n^{(+)}; A_{N_0}) = E_r(\triangle\gamma_n^{(+)} | A_{N_0}) P(A_{N_0}) = c_0(N+N_0)\int_0^\infty \frac{(\lambda y)^{N_0}}{N_0!} e^{-\lambda y} dH_r(y).$$

The cost structure is similar to the previous option $1 \leq s < N_0$. Namely, the costs consist of the storage costs, the replenishment costs, and the costs associated with the shortage of goods. Analytical forms of presentation of certain types of costs are also similar to the corresponding forms of the previous case. Consequently,

$$E_r(\triangle\gamma_n^{(-)} | A_{N_0}) = \sum_{i=1}^{N} \frac{c_1 i}{\lambda} + \underbrace{c_2(N+N_0)}_{\substack{\text{stock} \\ \text{replenishment} \\ \text{costs}}} + \underbrace{\sum_{i=1}^{N_0-1} \frac{c_3 i}{\lambda}}_{\substack{\text{costs related to the deficit} \\ \text{before the arrival of} \\ \text{the last order}}} + \underbrace{E_r^{(3*)}(\triangle\gamma_n^{(-)} | A_{N_0})}_{\substack{\text{costs related to the deficit} \\ \text{after the arrival of} \\ \text{the last order}}} =$$

$$= \frac{c_1 N(N+1)}{2\lambda} + c_2(N+N_0) + \frac{c_3 N_0(N_0-1)}{2\lambda} + E_r^{(3*)}(\triangle\gamma_n^{(-)} | A_{N_0}).$$

We pass from the conditional mathematical expectation to the mathematical expectation by the joint distribution

$$E_r(\triangle\gamma_n^{(-)} ; A_{N_0}) = E_r(\triangle\gamma_n^{(-)} | A_{N_0})P(A_{N_0}) = \left[\frac{c_1 N(N+1)}{2\lambda} + c_2(N+N_0) + \right.$$

$$\left. + \frac{c_3 N_0(N_0-1)}{2\lambda}\right]P(A_{N_0}) + E_r^{(3*)}(\triangle\gamma_n^{(-)} ; A_{N_0}).$$

In the standard way, we find an explicit expression for $E_r^{(3*)}(\triangle\gamma_n^{(-)} ; A_{N_0})$

$$E_r^{(3*)}(\triangle\gamma_n^{(-)} ; A_{N_0}) = c_3 N_0 \tau_{r,N_0} = c_3 N_0 \int_0^{\infty} \frac{\lambda^{N_0}}{N_0!} \int_x^{\infty} (z-x)^{N_0} e^{-\lambda z} dH_r(z) dx.$$

Thus,

$$E_r(\triangle\gamma_n^{(-)} ; A_{N_0}) = \left[\frac{c_1 N(N+1)}{2\lambda} + c_2(N+N_0) + \frac{c_3 N_0(N_0-1)}{2\lambda}\right]\int_0^{\infty} \frac{(\lambda y)^{N_0}}{N_0!} e^{-\lambda y} dH_r(y) +$$

$$+ c_3 N_0 \int_0^{\infty} \frac{\lambda^{N_0}}{N_0!} \int_x^{\infty} (z-x)^{N_0} e^{-\lambda z} dH_r(z) dx.$$

Based on the previous results, we obtain an analytical representation for the mathematical expectation of profit by the joint distribution with event $A_{N0}$:

$$E_r(\triangle\gamma_n; A_{N_0}) = E_r(\triangle\gamma_n^{(+)}; A_{N_0}) - E_r(\triangle\gamma_n^{(-)}; A_{N_0}) = c_0(N+N_0)\int_0^\infty \frac{(\lambda y)^{N_0}}{N_0!}e^{-\lambda y}dH_r(y) -$$

$$- \left[\frac{c_1 N(N+1)}{2\lambda} + c_2(N+N_0) + \frac{c_3 N_0(N_0-1)}{2\lambda}\right]\int_0^\infty \frac{(\lambda y)^{N_0}}{N_0!}e^{-\lambda y}dH_r(y) -$$

$$-c_3 N_0 \int_0^\infty \frac{\lambda^{N_0}}{N_0!}\int_x^\infty (z-x)^{N_0} e^{-\lambda z}dH_r(z)dx = \left[c_0(N+N_0) - \frac{c_1 N(N+1)}{2\lambda} - c_2(N+N_0) -\right.$$

$$\left.- \frac{c_3 N_0(N_0-1)}{2\lambda}\right]\int_0^\infty \frac{(\lambda y)^{N_0}}{N_0!}e^{-\lambda y}dH_r(y) - c_3 N_0 \int_0^\infty \frac{\lambda^{N_0}}{N_0!}\int_x^\infty (z-x)^{N_0}e^{-\lambda z}dH_r(z)dx. \quad (47)$$

8.3.4 Consider the case when $s > N_0$.

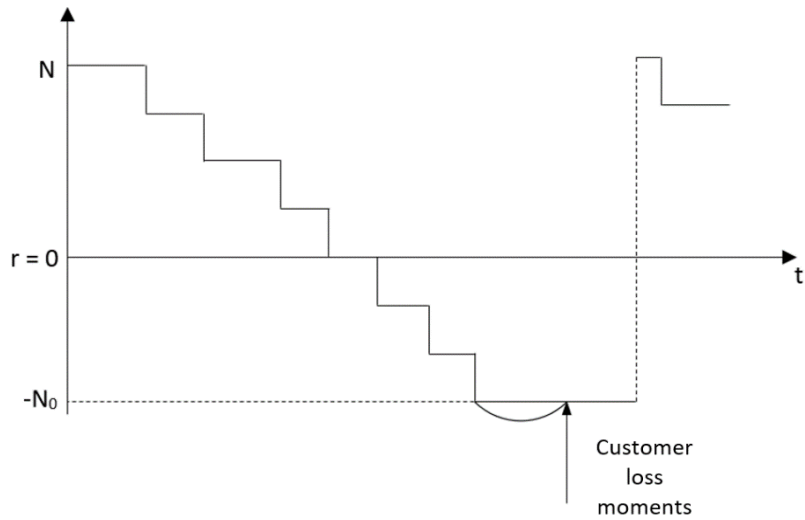

Fig. 18. Graphic illustration for the case 8.3.4

In this case, for any fixed $s > N_0$, the system will satisfy $N + N_0$ requirements for the goods, including N due to real stock and $N_0$ in the form of deferred demand. From here we find the formula for the mathematical expectation of income by the joint distribution

$$E_r(\triangle\gamma_n^{(+)}; A_s) = E_r(\triangle\gamma_n^{(+)} | A_s)P(A_s) = c_0(N+N_0)\int_0^\infty \frac{(\lambda y)^s}{s!}e^{-\lambda y}dH_r(y).$$

In the considered option, the cost structure is the most complex and contains all types of costs provided for in this mathematical model. We will separately comment on the form of costs associated with the deficit. These costs are formed from the following components:

1) Costs (penalties) that are formed during the period from the moment of deficit beginning to the moment when the deficit reaches the maximum allowable amount $N_0$. In other words, at the indicated moment, a demand with the number $N_0$ is received in the system, counting from the beginning of the delay period. These penalties accumulate linearly in time and depend on the volume of deficit.

2) Costs (penalties) that are formed during the period of time from the moment the maximum deficit is reached until the last requirement. Provided that the event $A_s$, s > N0 occurs, exactly s - $N_0$ requirements will be received for the specified period of time. The mathematical expectation of the duration of this period is equal to $\frac{s - N_0}{\lambda}$. Throughout this period, the deficit is $N_0$.

3) Costs (penalties) that are formed during the period of time from the moment of receipt of the last request to the end of the regeneration period. The conditional mathematical expectation of these costs is denoted by $E_r^{(3*)}(\triangle \gamma_n^{(-)} | A_s)$.

We note in addition that, provided that the event $A_s$, s > $N_0$, occurs during this regeneration period, s - $N_0$ requirements (which cannot be satisfied either as the real stock or in the form of deferred demand) will be lost. The costs associated with the irretrievable loss of customers are equal to $c_4^{(s-N_0)}$.

Таким образом, условное математическое ожидание расходов при условии осуществления события $A_s$, s > $N_0$, может быть представлено в виде

Thus, the conditional expectation of expenses under the condition of the event $A_s$, s > $N_0$, can be represented as

$$E_r(\triangle \gamma_n^{(-)} | A_s) = \underbrace{\sum_{i=1}^{N} \frac{c_1 i}{\lambda}}_{\substack{product \\ storage \\ costs}} + \underbrace{c_2(N + N_0)}_{\substack{stock \\ replenishment \\ costs}} + \underbrace{\sum_{i=1}^{N_0 - 1} \frac{c_3 i}{\lambda}}_{\substack{costs\ related\ to \\ the\ deficit\ before \\ order\ N_0}} + \underbrace{\frac{c_3 N_0}{\lambda}(s - N_0)}_{\substack{costs\ related\ to\ the\ deficit \\ after\ order\ N_0\ but \\ before\ the\ last\ order}} +$$

$$+ \underbrace{E_r^{(3*)}(\triangle \gamma_n^{(-)} | A_s)}_{\substack{costs\ related\ to\ the\ deficit \\ after\ the\ arrival \\ of\ the\ last\ order}} + \underbrace{c_4^{(s-N_0)}}_{\substack{costs\ of\ the\ loss \\ of\ clients}} = \frac{c_1 N(N+1)}{2\lambda} + c_2(N + N_0) + \frac{c_3 N_0(N_0 - 1)}{2\lambda} +$$

$$+ \frac{c_3 N_0}{\lambda}(s - N_0) + E_r^{(3*)}(\triangle \gamma_n^{(-)} | A_s) + c_4^{(s-N_0)}.$$

To obtain a representation of the corresponding mathematical expectation by the joint distribution with the event $A_s$, we multiply both sides of the last equality by the probability $P(A_s)$

$$E_r(\Delta\gamma_n^{(-)};A_s) = \left[\frac{c_1 N(N+1)}{2\lambda} + c_2(N+N_0) + \frac{c_3 N_0(2s-N_0-1)}{2\lambda} + c_4^{(s-N_0)}\right]P(A_s) +$$
$$+ E_r^{(3*)}(\Delta\gamma_n^{(-)};A_s).$$

We use the equality $P(A_s) = \int_0^\infty \frac{(\lambda y)^s}{s!} e^{-\lambda y} dH_r(y)$ and formula (23) for $E_r^{(3*)}(\Delta\gamma_n^{(-)};A_s)$:

$$E_r^{(3*)}(\Delta\gamma_n^{(-)};A_s) = c_3 N_0 \tau_{r,s} = c_3 N_0 \int_0^\infty \frac{\lambda^s}{s!} \int_x^\infty (z-x)^s e^{-\lambda z} dH_r(z) dx,$$

Then

$$E_r(\Delta\gamma_n^{(-)};A_s) = \left[\frac{c_1 N(N+1)}{2\lambda} + c_2(N+N_0) + \frac{c_3 N_0(2s-N_0-1)}{2\lambda} + \right.$$
$$\left. + c_4^{(s-N_0)}\right] \int_0^\infty \frac{(\lambda y)^s}{s!} e^{-\lambda y} dH_r(y) + c_3 N_0 \int_0^\infty \frac{\lambda^s}{s!} \int_x^\infty (z-x)^s e^{-\lambda z} dH_r(z) dx.$$

Thus, the mathematical expectation of profit for the period by the joint distribution with the event $A_s$ can be written in the following form.

$$E_r(\Delta\gamma_n;A_s) = E_r(\Delta\gamma_n^{(+)};A_s) - E_r(\Delta\gamma_n^{(-)};A_s) = c_0(N+N_0)\int_0^\infty \frac{(\lambda y)^s}{s!} e^{-\lambda y} dH_r(y) -$$
$$-\left[\frac{c_1 N(N+1)}{2\lambda} + c_2(N+N_0) + \frac{c_3 N_0(2s-N_0-1)}{2\lambda} + c_4^{(s-N_0)}\right]\int_0^\infty \frac{(\lambda y)^s}{s!} e^{-\lambda y} dH_r(y) -$$
$$- c_3 N_0 \int_0^\infty \frac{\lambda^s}{s!} \int_x^\infty (z-x)^s e^{-\lambda z} dH_r(z) dx = \left[c_0(N+N_0) - \frac{c_1 N(N+1)}{2\lambda} - c_2(N+N_0) - \right.$$
$$\left. - \frac{c_3 N_0(2s-N_0-1)}{2\lambda} - c_4^{(s-N_0)}\right]\int_0^\infty \frac{(\lambda y)^s}{s!} e^{-\lambda y} dH_r(y) - c_3 N_0 \int_0^\infty \frac{\lambda^s}{s!} \int_x^\infty (z-x)^s e^{-\lambda z} dH_r(z) dx. \quad (48)$$

Based on relation (27), taking into account equalities (45) - (48), we obtain an explicit representation for the mathematical expectation of profit, provided that on the regeneration period the control parameter takes a fixed value r = 0. As before, we underline the parts of the right side of the next relation related to the corresponding values of the parameter s.

$$E_r(\Delta \gamma_n) = \underbrace{\left[c_0 N - \frac{c_1 N(N+1)}{2\lambda} - c_2 N\right]\int_0^\infty e^{-\lambda y} dH_r(y)}_{for\ s=0} + \underbrace{\sum_{i=1}^{N_0-1}\left[\left[c_0(N+i) - \frac{c_1 N(N+1)}{2\lambda}\right.\right.}_{}$$

$$\underbrace{\left.\left.- c_2(N+i) - \frac{c_3 i(i-1)}{2\lambda}\right]\int_0^\infty \frac{(\lambda y)^i}{i!} e^{-\lambda y} dH_r(y) - c_3 i \int_0^\infty \frac{\lambda^i}{i!}\int_x^\infty (z-x)^i e^{-\lambda z} dH_r(z) dx\right]}_{for\ 1\leq s < N_0} +$$

$$+ \underbrace{\left[c_0(N+N_0) - \frac{c_1 N(N+1)}{2\lambda} - c_2(N+N_0) - \frac{c_3 N_0(N_0-1)}{2\lambda}\right]\int_0^\infty \frac{(\lambda y)^{N_0}}{N_0!} e^{-\lambda y} dH_r(y)}_{for\ s=N_0}$$

$$\underbrace{- c_3 N_0 \int_0^\infty \frac{\lambda^{N_0}}{N_0!}\int_x^\infty (z-x)^{N_0} e^{-\lambda z} dH_r(z) dx}_{} + \underbrace{\sum_{i=N_0+1}^\infty \left[(c_0(N+N_0) - \frac{c_1 N(N+1)}{2\lambda} - c_2(N+N_0) -\right.}_{for\ s>N_0}$$

$$\underbrace{\left.- \frac{c_3 N_0(2i-N_0-1)}{2\lambda} - c_4^{(i-N_0)}\right]\int_0^\infty \frac{(\lambda y)^i}{i!} e^{-\lambda y} dH_r(y) - c_3 N_0 \int_0^\infty \frac{\lambda^i}{i!}\int_x^\infty (z-x)^i e^{-\lambda z} dH_r(z) dx\right].}_{}$$

We can now write an explicit representation for the control efficiency indicator $I_r$. We note again that the resulting representation is valid for $r = 0$.

$$I_r = \frac{1}{\frac{N}{\lambda} + \int_0^\infty (1-H_r(y))dy}\left[\left[c_0 N - \frac{c_1 N(N+1)}{2\lambda} - c_2 N\right]\int_0^\infty e^{-\lambda y} dH_r(y) + \sum_{i=1}^{N_0-1}\left[\left[c_0(N+i) -\right.\right.\right.$$

$$\left.\frac{c_1 N(N+1)}{2\lambda} - c_2(N+i) - \frac{c_3 i(i-1)}{2\lambda}\right]\int_0^\infty \frac{(\lambda y)^i}{i!} e^{-\lambda y} dH_r(y) - c_3 i\int_0^\infty \frac{\lambda^i}{i!}\int_x^\infty (z-x)^i e^{-\lambda z} dH_r(z) dx\right] +$$

$$+ \left[c_0(N+N_0) - \frac{c_1 N(N+1)}{2\lambda} - c_2(N+N_0) - \frac{c_3 N_0(N_0-1)}{2\lambda}\right]\int_0^\infty \frac{(\lambda y)^{N_0}}{N_0!} e^{-\lambda y} dH_r(y) -$$

$$- c_3 N_0 \int_0^\infty \frac{\lambda^{N_0}}{N_0!}\int_x^\infty (z-x)^{N_0} e^{-\lambda z} dH_r(z) dx + \sum_{i=N_0+1}^\infty \left[\left[c_0(N+N_0) - \frac{c_1 N(N+1)}{2\lambda} - c_2(N+N_0) -\right.\right.$$

$$\left.\left.\left.- \frac{c_3 N_0(2i-N_0-1)}{2\lambda} - c_4^{(i-N_0)}\right]\int_0^\infty \frac{(\lambda y)^i}{i!} e^{-\lambda y} dH_r(y) - c_3 N_0 \int_0^\infty \frac{\lambda^i}{i!}\int_x^\infty (z-x)^i e^{-\lambda z} dH_r(z) dx\right]\right].(49)$$

## 8.4 The variant when $-N_0 < r \leq -1$

### 8.4.1 Consider the case when s = 0.

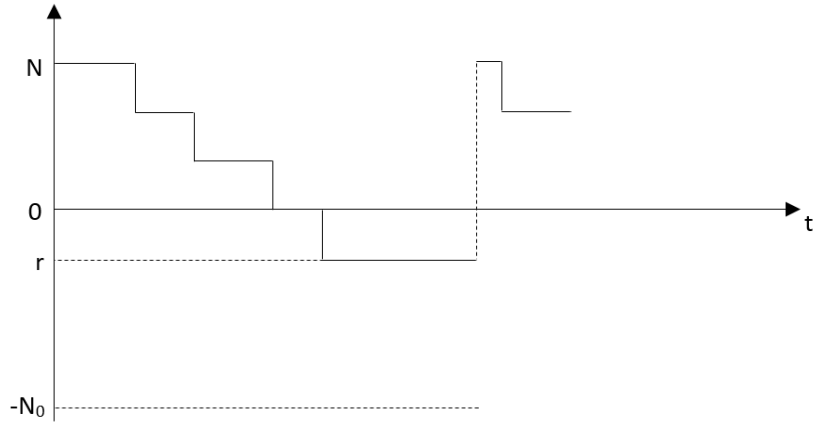

Fig. 19. Graphic illustration for the case 8.4.1

Since in this case, when event $A_0$ occurs, no requirements were received by the system during the delay period, N – r units of goods will be consumed over the entire regeneration interval, and N – r > N. Of this number, N units will be delivered due to the available real stock, and (-r) units according to the deferred demand scheme. Thus, the mathematical expectation of income by the joint distribution with event $A_0$ is determined by the relation

$$E_r(\triangle\gamma_n^{(+)}; A_0) = E_r(\triangle\gamma_n^{(+)} | A_0)P(A_0) = c_0(N-r)\int_0^\infty e^{-\lambda y}dH_r(y).$$

The costs for the considered period of regeneration are the sum of the storage costs, the costs necessary to replenish the stock, and the costs associated with the deficit. Moreover, the costs associated with the deficit consist of costs incurred in the time interval from the start of the regeneration period to the time of the replenishment order (or the start of the delay period) and the costs incurred in the delivery delay period. In the first of the indicated time intervals, the deficit changes from a value of 1 to a value of (-r - 1), and in the second it keeps a fixed value (-r). Note that the deficit is conveniently calculated in positive units. From this we obtain the following representation for the conditional mathematical expectation of costs during the regeneration period

$$E_r(\triangle\gamma_n^{(-)} \mid A_0) = \underbrace{\sum_{i=1}^{N}\frac{c_1 i}{\lambda}}_{\substack{product\\storage\\costs}} + \underbrace{c_2(N-r)}_{\substack{stock\\replenishment\\costs}} + \underbrace{\sum_{i=1}^{-r-1}\frac{c_3 i}{\lambda}}_{\substack{costs\ related\ to\ the\ deficit\\before\ the\ beggining\ of\\the\ delay\ period}} + \underbrace{c_3(-r)\int_0^{\infty}(1-H_r(y))dy}_{\substack{costs\ related\ to\ the\ deficit\\after\ the\ beggining\ of\\the\ delay\ period}}.$$

Therefore, the mathematical expectation of expenses by joint distribution with event $A_0$ is determined by the formula

$$E_r(\triangle\gamma_n^{(-)}; A_0) = E_r(\triangle\gamma_n^{(-)} \mid A_0)P(A_0) = \left[\frac{c_1 N(N+1)}{2\lambda} + c_2(N-r) + \frac{c_3 r(r+1)}{2\lambda} + \right.$$

$$\left. + c_3(-r)\int_0^{\infty}(1-H_r(y))dy\right]\int_0^{\infty}e^{-\lambda y}dH_r(y).$$

Thus, the corresponding mathematical expectation of profit during the regeneration period has the form

$$E_r(\triangle\gamma_n; A_0) = E_r(\triangle\gamma_n^{(+)}; A_0) - E_r(\triangle\gamma_n^{(-)}; A_0) = c_0(N-r)\int_0^{\infty}e^{-\lambda y}dH_r(y) - \left[\frac{c_1 N(N+1)}{2\lambda} + \right.$$

$$\left. + c_2(N-r) + \frac{c_3 r(r+1)}{2\lambda} + c_3(-r)\int_0^{\infty}(1-H_r(y))dy\right]\int_0^{\infty}e^{-\lambda y}dH_r(y) = \left[c_0(N-r) - \frac{c_1 N(N+1)}{2\lambda} - \right.$$

$$\left. - c_2(N-r) - \frac{c_3 r(r+1)}{2\lambda} + c_3 r\int_0^{\infty}(1-H_r(y))dy\right]\int_0^{\infty}e^{-\lambda y}dH_r(y). \qquad (50)$$

8.4.2  Consider the case when $1 \leq s < N_0 + r$.

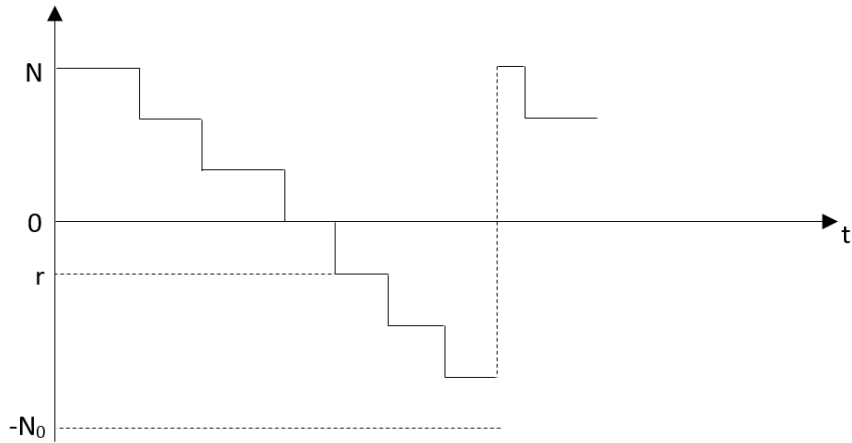

Fig. 20. Graphic illustration for the case 8.4.2

In this case, as in the previous one, by the time of the order or the start of the delivery delay, a deficit in the amount of (-r) is formed in the system. When $s \geq 1$ requirements are received during the delay period, the deficit reaches s - r. At the same time, it is assumed that the resulting total deficit does not reach its maximum permissible level $N_0$, that is, all incoming requirements will be satisfied due to real stock and deferred demand. Thus, N + (s - r) of the received requirements will be paid, and the mathematical expectation of income by the joint distribution with the event $A_s$ is determined by the following formula

$$E_r(\Delta\gamma_n^{(+)}; A_s) = E_r(\Delta\gamma_n^{(+)} | A_s)P(A_s) = c_0(N - r + s)\int_0^\infty \frac{(\lambda y)^s}{s!} e^{-\lambda y} dH_r(y).$$

In this case, the costs consist of the storage costs, the costs necessary to replenish the stock, and the costs due to the deficit. The costs due to the deficit can be represented as the sum of the costs in the time interval from the moment the deficit was formed until the last request arrived and the costs in the time interval from the time the last request arrived until the end of the regeneration period. Note that to replenish the stock, N - r + s units of goods are needed, and the deficit at the time of receipt of the last demand reaches the level of (s - r). From here we get a representation for the conditional expectation of expenses

$$E_r(\Delta\gamma_n^{(-)} | A_s) = \underbrace{\sum_{i=1}^{N} \frac{c_1 i}{\lambda}}_{\substack{product \\ storage \\ costs}} + \underbrace{c_2(N - r + s)}_{\substack{stock \\ replenishment \\ costs}} + \underbrace{\sum_{i=1}^{-r+s-1} \frac{c_3 i}{\lambda}}_{\substack{costs\ related\ to\ the\ deficit \\ before\ the\ arrival\ of \\ the\ last\ order}} + \underbrace{E_r^{(3*)}(\Delta\gamma_n^{(-)} | A_s)}_{\substack{costs\ related\ to\ the\ deficit \\ after\ the\ arrival\ of \\ the\ last\ order}}.$$

Passing to the mathematical expectation of expenses during the regeneration period by the joint distribution with event $A_s$, and considering that

$$\sum_{i=1}^{N} \frac{c_1 i}{\lambda} = \frac{c_1 N(N+1)}{2\lambda} \quad \text{and} \quad \sum_{i=1}^{-r+s-1} \frac{c_3 i}{\lambda} = \frac{c_3(-r+s)(-r+s-1)}{2\lambda},$$

we obtain

$$E_r(\Delta\gamma_n^{(-)}; A_s) = E_r(\Delta\gamma_n^{(-)} | A_s)P(A_s) = \left[\frac{c_1 N(N+1)}{2\lambda} + c_2(N - r + s) + \frac{c_3(-r+s)(-r+s-1)}{2\lambda}\right]P(A_s) + E_r^{(3*)}(\Delta\gamma_n^{(-)}; A_s).$$

To find an explicit representation of the mathematical expectation of costs due to the deficit and arising in the time interval from the moment the last demand arrived until the end of the regeneration period, we use the statement of Theorem 3. Note that this theorem will be used in the future to find similar cost characteristics associated with the deficit

$$E_r^{(3*)}(\triangle \gamma_n^{(-)}; A_s) = c_3(-r+s)\tau_{r,s} = c_3(-r+s)\int_0^\infty \frac{\lambda^s}{s!}\int_x^\infty (z-x)^s e^{-\lambda z} dH_r(z) dx.$$

This implies

$$E_r(\triangle \gamma_n^{(-)}; A_s) = \left[\frac{c_1 N(N+1)}{2\lambda} + c_2(N-r+s) + \frac{c_3(-r+s)(-r+s-1)}{2\lambda}\right]\int_0^\infty \frac{(\lambda y)^s}{s!} e^{-\lambda y} dH_r(y) +$$

$$+ c_3(-r+s)\int_0^\infty \frac{\lambda^s}{s!}\int_x^\infty (z-x)^s e^{-\lambda z} dH_r(z) dx.$$

Thus, the corresponding mathematical expectation of profit for the case in question will take the form:

$$E_r(\triangle \gamma_n; A_s) = E_r(\triangle \gamma_n^{(+)}; A_s) - E_r(\triangle \gamma_n^{(-)}; A_s) = c_0(N-r+s)\int_0^\infty \frac{(\lambda y)^s}{s!} e^{-\lambda y} dH_r(y) -$$

$$-\left[\frac{c_1 N(N+1)}{2\lambda} + c_2(N-r+s) + \frac{c_3(-r+s)(-r+s-1)}{2\lambda}\right]\int_0^\infty \frac{(\lambda y)^s}{s!} e^{-\lambda y} dH_r(y) -$$

$$-c_3(-r+s)\int_0^\infty \frac{\lambda^s}{s!}\int_x^\infty (z-x)^s e^{-\lambda z} dH_r(z)dx = \left[c_0(N-r+s) - \frac{c_1 N(N+1)}{2\lambda} - c_2(N-r+s) - \right.$$

$$\left. - \frac{c_3(-r+s)(-r+s-1)}{2\lambda}\right]\int_0^\infty \frac{(\lambda y)^s}{s!} e^{-\lambda y} dH_r(y) - c_3(-r+s)\int_0^\infty \frac{\lambda^s}{s!}\int_x^\infty (z-x)^s e^{-\lambda z} dH_r(z) dx. \quad (51)$$

8.4.3 Consider the case when $s = N_0 + r$.

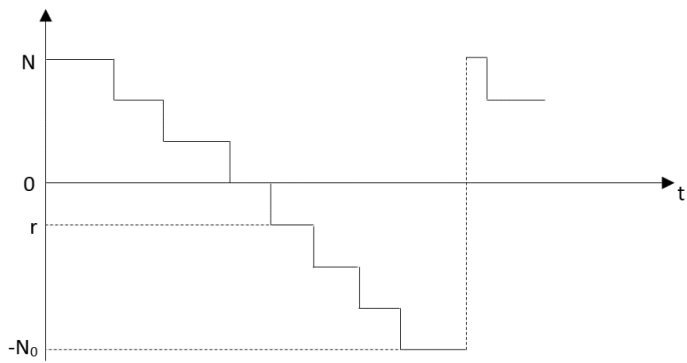

Fig. 21. Graphic illustration for the case 8.4.3

In this case, during the regeneration period, the system will receive exactly $N + N_0$ requirements, of which $N - r$ will arrive in the time interval before the start of the delay period (including the moment the beginning of this period), and the remaining $s = N_0$

+ r during the delivery delay period. All incoming requirements will be satisfied using real stock and according to the deferred demand scheme. Thus,

$$E_r(\triangle\gamma_n^{(+)}; A_{N_0+r}) = E_r(\triangle\gamma_n^{(+)} | A_{N_0+r})P(A_{N_0+r}) = c_0(N+N_0)\int_0^\infty \frac{(\lambda y)^{N_0+r}}{(N_0+r)!}e^{-\lambda y}dH_r(y).$$

The general structure of the conditional mathematical expectation of costs is completely similar to the previous case. We only note that in this case, to replenish the stock, $N + N_0$ units of goods are needed, and the deficit at the moment of receipt of the last demand reaches the level of $N_0$ and remains at this level until the end of the regeneration period. From here we get

$$E_r(\triangle\gamma_n^{(-)} | A_{N_0+r}) = \underbrace{\sum_{i=1}^{N}\frac{c_1 i}{\lambda}}_{\substack{\text{product}\\\text{storage}\\\text{costs}}} + \underbrace{c_2(N+N_0)}_{\substack{\text{stock}\\\text{replenishment}\\\text{costs}}} + \underbrace{\sum_{i=1}^{N_0-1}\frac{c_3 i}{\lambda}}_{\substack{\text{costs related to the deficit}\\\text{before the arrival of}\\\text{the last order}}} + \underbrace{E_r^{(3*)}(\triangle\gamma_n^{(-)} | A_{N_0+r})}_{\substack{\text{costs related to the deficit}\\\text{after the arrival of}\\\text{the last order}}}.$$

The last relation directly implies the representation for the mathematical expectation of the costs by the joint distribution with event $A_{N0+r}$

$$E_r(\triangle\gamma_n^{(-)}; A_{N_0+r}) = E_r(\triangle\gamma_n^{(-)} | A_{N_0+r})P(A_{N_0+r}) = \left[\frac{c_1 N(N+1)}{2\lambda} + c_2(N+N_0) + \right.$$

$$\left. + \frac{c_3 N_0(N_0-1)}{2\lambda}\right]\int_0^\infty \frac{(\lambda y)^{N_0+r}}{(N_0+r)!}e^{-\lambda y}dH_r(y) + E_r^{(3*)}(\triangle\gamma_n^{(-)}; A_{N_0+r}).$$

To calculate the mathematical expectation $E_r^{(3*)}(\triangle\gamma_n^{(-)}; A_{N_0+r})$, we again use the statement of Theorem 3

$$E_r^{(3*)}(\triangle\gamma_n^{(-)}; A_{N_0+r}) = c_3 N_0 \tau_{r,N_0+r} = c_3 N_0 \int_0^\infty \frac{\lambda^{N_0+r}}{(N_0+r)!}\int_x^\infty (z-x)^{N_0+r}e^{-\lambda z}dH_r(z)dx.$$

Thus, the expression for $E_r(\triangle\gamma_n^{(-)}; A_{N_0+r})$ can be rewritten as:

$$E_r(\triangle\gamma_n^{(-)}; A_{N_0+r}) = E_r(\triangle\gamma_n^{(-)} | A_{N_0+r})P(A_{N_0+r}) = (\frac{c_1 N(N+1)}{2\lambda} + c_2(N+N_0) +$$

$$+ \frac{c_3 N_0(N_0-1)}{2\lambda})\int_0^\infty \frac{(\lambda y)^{N_0+r}}{(N_0+r)!}e^{-\lambda y}dH_r(y) + c_3 N_0 \int_0^\infty \frac{\lambda^{N_0+r}}{(N_0+r)!}\int_x^\infty (z-x)^{N_0+r}e^{-\lambda z}dH_r(z)dx.$$

Using the above relations, we write out the final representation for the mathematical expectation of profit during the regeneration period by the joint distribution with event $A_{N_0+r}$:

$$E_r(\triangle\gamma_n; A_{N_0+r}) = E_r(\triangle\gamma_n^{(+)}; A_{N_0+r}) - E_r(\triangle\gamma_n^{(-)}; A_{N_0+r}) = c_0(N+N_0)\int_0^\infty \frac{(\lambda y)^{N_0+r}}{(N_0+r)!}e^{-\lambda y}dH_r(y) -$$

$$-\left[\frac{c_1 N(N+1)}{2\lambda} + c_2(N+N_0) + \frac{c_3 N_0(N_0-1)}{2\lambda}\right]\int_0^\infty \frac{(\lambda y)^{N_0+r}}{(N_0+r)!}e^{-\lambda y}dH_r(y) -$$

$$-c_3 N_0 \int_0^\infty \frac{\lambda^{N_0+r}}{(N_0+r)!}\int_x^\infty (z-x)^{N_0+r}e^{-\lambda z}dH_r(z)dx = \left[c_0(N+N_0) - \frac{c_1 N(N+1)}{2\lambda} - c_2(N+N_0) - \right.$$

$$\left. - \frac{c_3 N_0(N_0-1)}{2\lambda}\right]\int_0^\infty \frac{(\lambda y)^{N_0+r}}{(N_0+r)!}e^{-\lambda y}dH_r(y) - c_3 N_0 \int_0^\infty \frac{\lambda^{N_0+r}}{(N_0+r)!}\int_x^\infty (z-x)^{N_0+r}e^{-\lambda z}dH_r(z)dx. \quad (52)$$

8.4.4 Consider the case when $s > N_0 + r$.

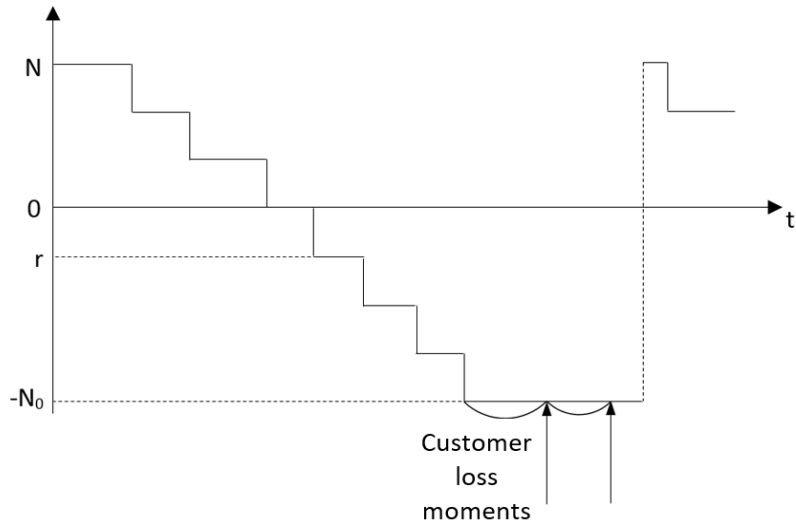

Fig. 22. Graphic illustration for the case 8.4.4

Under the condition of the realization of event $A_s$, $s > N_0 + r$, the number of requirements received by the system during the regeneration period will exceed $N + N_0$. According to the model, $N + N_0$ requirements of them will be satisfied, and the mathematical expectation of income is determined by the formula

$$E_r(\triangle\gamma_n^{(+)}; A_s) = E_r(\triangle\gamma_n^{(+)} | A_s)P(A_s) = c_0(N+N_0)\int_0^\infty \frac{(\lambda y)^s}{s!}e^{-\lambda y}dH_r(y).$$

Costs in this case have the most complex structure and consist of four main components:

1) real stock storage costs;
2) stock replenishment costs;
3) costs related to product deficit;
4) costs resulting from loss of clients (orders).

At the same time, the costs related to the product deficit should be divided into three parts with different analytical expressions:

3.1) The costs incurred in the time interval from the moment of deficit formation to the moment when the deficit reaches the level of $N_0$. Note that at this moment the request with the number $N_0 + r$ enters the system, counting from the moment the delay period begins.

3.2) The costs incurred in the time interval from the moment of receipt of the request with the number $N_0 + r$ to the moment of receipt of the last request on the regeneration period. Throughout this interval, the deficit remains at the level of $N_0$. The distribution of this interval is Erlangian order $s - r - N_0$, since the last incoming request has the number s, counting from the beginning of the delay period.

3.3) The costs incurred in the time interval from the moment the last request arrives until the end of the regeneration period. The conditional mathematical expectation of these costs, as before, is denoted by $E_r^{(3*)}(\triangle\gamma_n^{(-)} | A_s)$.

We note in addition that in the case under consideration, the first $N_0 + r$ requirements received by the system from the moment beginning of the delay period will be satisfied, and the remaining $s - r - N_0$ will be lost. Thus, the costs associated with the loss of customers, are equal to $c_4^{(s-r-N_0)}$.

Based on the above comments, we obtain a representation for the conditional mathematical expectation of costs

$$E_r(\triangle\gamma_n^{(-)} | A_s) = \underbrace{\sum_{i=1}^{N}\frac{c_1 i}{\lambda}}_{\substack{product \\ storage \\ costs}} + \underbrace{c_2(N + N_0)}_{\substack{stock \\ replenishment \\ costs}} + \underbrace{\sum_{i=1}^{N_0-1}\frac{c_3 i}{\lambda}}_{\substack{costs\ related\ to \\ the\ deficit\ before \\ order\ r+N_0}} + \underbrace{\frac{c_3 N_0}{\lambda}(s - r - N_0)}_{\substack{costs\ related\ to\ the\ deficit \\ after\ order\ r+N_0\ but \\ before\ the\ last\ order}} +$$

$$+ \underbrace{E_r^{(3*)}(\triangle\gamma_n^{(-)} | A_s)}_{\substack{costs\ related\ to\ the\ deficit \\ after\ the\ arrival \\ of\ the\ last\ order}} + \underbrace{c_4^{(s-r-N_0)}}_{\substack{costs\ of\ the\ loss \\ of\ clients}}.$$

This implies

$$E_r(\Delta\gamma_n^{(-)};A_s) = E_r(\Delta\gamma_n^{(-)} | A_s)P(A_s) = \left[\frac{c_1 N(N+1)}{2\lambda} + c_2(N+N_0) + \frac{c_3 N_0(N_0-1)}{2\lambda} + \right.$$

$$\left. + \frac{c_3 N_0}{\lambda}(s-r-N_0) + c_4^{(s-r-N_0)}\right]\int_0^\infty \frac{(\lambda y)^s}{s!} e^{-\lambda y} dH_r(y) + E_r^{(3*)}(\Delta\gamma_n^{(-)};A_s).$$

The analytical representation of the mathematical expectation $E_r^{(3*)}(\Delta\gamma_n^{(-)};A_s)$, as in the previous cases, is determined on the basis of the statement of Theorem 3 (formula (23)):

$$E_r(\Delta\gamma_n^{(-)};A_s) = E_r(\Delta\gamma_n^{(-)} | A_s)P(A_s) = (\frac{c_1 N(N+1)}{2\lambda} + c_2(N+N_0) + \frac{c_3 N_0(N_0-1)}{2\lambda} +$$

$$+ \frac{c_3 N_0}{\lambda}(s-r-N_0) + c_4^{(s-r-N_0)})\int_0^\infty \frac{(\lambda y)^s}{s!} e^{-\lambda y} dH_r(y) + c_3 N_0 \int_0^\infty \frac{\lambda^s}{s!}\int_x^\infty (z-x)^s e^{-\lambda z} dH_r(z)dx.$$

Let us summarize the analysis of this case and write out the expression for $E_r(\Delta\gamma_n;A_s)$:

$$E_r(\Delta\gamma_n;A_s) = E_r(\Delta\gamma_n^{(+)};A_s) - E_r(\Delta\gamma_n^{(-)};A_s) = c_0(N+N_0)\int_0^\infty \frac{(\lambda y)^s}{s!} e^{-\lambda y} dH_r(y) - \left[\frac{c_1 N(N+1)}{2\lambda} + \right.$$

$$+ c_2(N+N_0) + \frac{c_3 N_0(2s-2r-N_0-1)}{2\lambda} + c_4^{(s-r-N_0)}\left.\right]\int_0^\infty \frac{(\lambda y)^s}{s!} e^{-\lambda y} dH_r(y) -$$

$$- c_3 N_0 \int_0^\infty \frac{\lambda^s}{s!}\int_x^\infty (z-x)^s e^{-\lambda z} dH_r(z)dx = \left[c_0(N+N_0) - \frac{c_1 N(N+1)}{2\lambda} - c_2(N+N_0) - \right.$$

$$\left. - \frac{c_3 N_0(2s-2r-N_0-1)}{2\lambda} - c_4^{(s-r-N_0)}\right]\int_0^\infty \frac{(\lambda y)^s}{s!} e^{-\lambda y} dH_r(y) - c_3 N_0 \int_0^\infty \frac{\lambda^s}{s!}\int_x^\infty (z-x)^s e^{-\lambda z} dH_r(z)dx. \quad (53)$$

From relation (27), taking into account formulas (50) - (53), we obtain an explicit representation for the mathematical expectation of profit during the regeneration period, provided that the decision r, $-N_0 < r \leq -1$ was made at this period. As before, for the convenience of perceiving such a cumbersome formula, we highlight in its parts corresponding to different ranges of the parameter s.

$$E_r(\Delta\gamma_n) = \underbrace{(c_0(N-r) - \frac{c_1 N(N+1)}{2\lambda} - c_2(N-r) - \frac{c_3 r(r+1)}{2\lambda} + c_3 r \int_0^\infty (1-H_r(y))dy)\int_0^\infty e^{-\lambda y}dH_r(y)}_{\text{for } s=0} +$$

$$+ \underbrace{\sum_{i=1}^{N_0+r-1}((c_0(N-r+i) - \frac{c_1 N(N+1)}{2\lambda} - c_2(N-r+i) - \frac{c_3(-r+i)(-r+i-1)}{2\lambda})\int_0^\infty \frac{(\lambda y)^i}{i!}e^{-\lambda y}dH_r(y)}_{\text{for } 1 \le s < N_0+r} -$$

$$-c_3(-r+i)\int_0^\infty \frac{\lambda^i}{i!}\int_x^\infty (z-x)^i e^{-\lambda z}dH_r(z)dx) + \underbrace{(c_0(N+N_0) - \frac{c_1 N(N+1)}{2\lambda} - c_2(N+N_0) -}_{\text{for } s=N_0+r}$$

$$-\frac{c_3 N_0(N_0-1)}{2\lambda})\int_0^\infty \frac{(\lambda y)^{N_0+r}}{(N_0+r)!}e^{-\lambda y}dH_r(y) - c_3 N_0 \int_0^\infty \frac{\lambda^{N_0+r}}{(N_0+r)!}\int_x^\infty (z-x)^{N_0+r}e^{-\lambda z}dH_r(z)dx +$$

$$+ \underbrace{\sum_{i=N_0+r+1}^\infty ((c_0(N+N_0) - \frac{c_1 N(N+1)}{2\lambda} - c_2(N+N_0) - \frac{c_3 N_0(2i-2r-N_0-1)}{2\lambda} -}_{\text{for } s > N_0+r}$$

$$-c_4^{(i-r-N_0)})\int_0^\infty \frac{(\lambda y)^i}{i!}e^{-\lambda y}dH_r(y) - c_3 N_0 \int_0^\infty \frac{\lambda^i}{i!}\int_x^\infty (z-x)^i e^{-\lambda z}dH_r(z)dx).$$

Now we can get an explicit representation for the control efficiency indicator $I_r$ for the variant when the control parameter r takes values in the region $-N_0 < r \le -1$. It should be noted that the mathematical expectation of the duration of the regeneration period is determined by the following relation: $E_r(\Delta t_n) = \frac{N-r}{\lambda} + \int_0^\infty (1-H_r(y))dy$.

$$I_r = \frac{1}{\frac{N-r}{\lambda} + \int_0^\infty (1-H_r(y))dy}\left\{\left[c_0(N-r) - \frac{c_1 N(N+1)}{2\lambda} - c_2(N-r) - \frac{c_3 r(r+1)}{2\lambda} + \right.\right.$$

$$\left.+ c_3 r\int_0^\infty (1-H_r(y))dy\right]\int_0^\infty e^{-\lambda y}dH_r(y) + \sum_{i=1}^{N_0+r-1}\left[\left[c_0(N-r+i) - \frac{c_1 N(N+1)}{2\lambda} -\right.\right.$$

$$\left.\left.- c_2(N-r+i) - \frac{c_3(-r+i)(-r+i-1)}{2\lambda}\right]\int_0^\infty \frac{(\lambda y)^i}{i!}e^{-\lambda y}dH_r(y) -\right.$$

$$\left.- c_3(-r+i)\int_0^\infty \frac{\lambda^i}{i!}\int_x^\infty (z-x)^i e^{-\lambda z}dH_r(z)dx\right] + \left[c_0(N+N_0) - \frac{c_1 N(N+1)}{2\lambda} - c_2(N+N_0) -\right.$$

$$\left.- \frac{c_3 N_0(N_0-1)}{2\lambda}\right]\int_0^\infty \frac{(\lambda y)^{N_0+r}}{(N_0+r)!}e^{-\lambda y}dH_r(y) - c_3 N_0\int_0^\infty \frac{\lambda^{N_0+r}}{(N_0+r)!}\int_x^\infty (z-x)^{N_0+r}e^{-\lambda z}dH_r(z)dx +$$

$$+ \sum_{i=N_0+r+1}^{\infty} \left[ \left[ c_0(N+N_0) - \frac{c_1 N(N+1)}{2\lambda} - c_2(N+N_0) - \frac{c_3 N_0(2i-2r-N_0-1)}{2\lambda} - \right.\right.$$
$$\left.\left. -c_4^{(i-r-N_0)} \right] \int_0^{\infty} \frac{(\lambda y)^i}{i!} e^{-\lambda y} dH_r(y) - c_3 N_0 \int_0^{\infty} \frac{\lambda^i}{i!} \int_x^{\infty} (z-x)^i e^{-\lambda z} dH_r(z) dx \right]. \quad (54)$$

## 8.5 The variant when r = -N₀

### 8.5.1 Consider the case when s = 0.

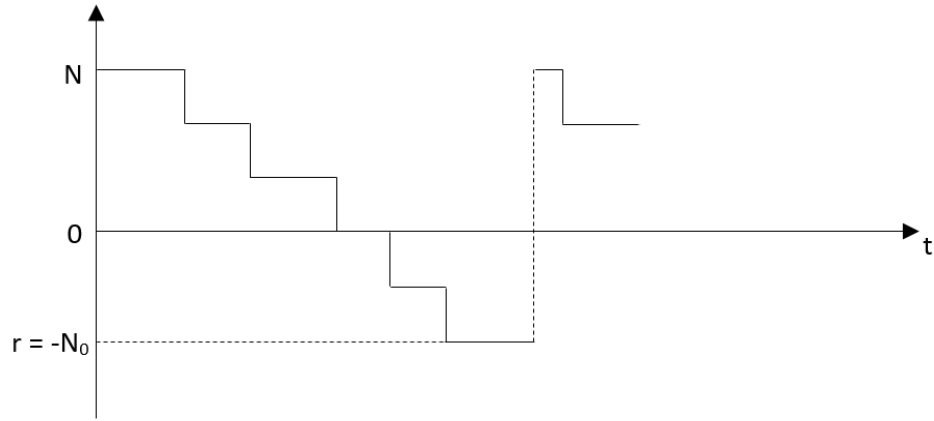

Fig. 23. Graphic illustration for the case 8.5.1

We consider separately the variant when the control parameter r takes the value -$N_0$. Within this variant, a replenishment order will be placed when the inventory level in the system reaches the minimum acceptable value (-$N_0$). At this point, the period of delay in delivery will begin. Note that then it is advisable to consider two mutually exclusive cases related to the receipt of customers during the delivery delay period:

1) Not a single client arrives, that is, event $A_0$ occurs;
2) A certain number of requests arrive, that is, one event occurs from the system of incompatible events $A_s$, $s \geq 1$. At the same time, the incoming requirements will be lost.

If event $A_0$ is realized, then all requirements received by the system during the regeneration period will be satisfied due to real stock and deferred demand. The number of these requirements is N + $N_0$. From here $E_r(\Delta \gamma_n^{(+)} | A_0) = c_0(N+N_0)$, and further, we obtain the representation for the mathematical expectation of income by the joint distribution with event $A_0$:

$$E_r(\triangle\gamma_n^{(+)}; A_0) = E_r(\triangle\gamma_n^{(+)} \mid A_0)P(A_0) = c_0(N + N_0)\int_0^\infty e^{-\lambda y}dH_r(y).$$

The costs in this case are the sum of the product storage costs, the costs necessary to replenish the stock, and the costs associated with the deficit of goods. The costs associated with the deficit are the sum of the costs in the time interval from the moment the deficit is formed to the moment when the deficit reaches the maximum allowable level $N_0$, and the costs in the time interval from the moment the deficit reaches the level of $N_0$ until the end of the regeneration period. Note that the moment when the deficit reaches the level of $N_0$ coincides with the moment of ordering or the beginning of the delivery delay period. During the entire delay period, the deficit remains at the level of $N_0$. Based on these comments, we find

$$E_r(\triangle\gamma_n^{(-)} \mid A_0) = \underbrace{\sum_{i=1}^{N}\frac{c_1 i}{\lambda}}_{\substack{product \\ storage \\ costs}} + \underbrace{c_2(N+N_0)}_{\substack{stock \\ replenishment \\ costs}} + \underbrace{\sum_{i=1}^{N_0-1}\frac{c_3 i}{\lambda}}_{\substack{costs\,related\,to\,the\,deficit \\ before\,the\,beggining\,of \\ the\,delay\,period}} + \underbrace{c_3 N_0 \int_0^\infty (1 - H_r(y))dy}_{\substack{costs\,related\,to\,the\,deficit \\ after\,the\,beggining\,of \\ the\,delay\,period}}.$$

Given that

$$\sum_{i=1}^{N}\frac{c_1 i}{\lambda} = \frac{c_1 N(N+1)}{2\lambda} \quad\text{and}\quad \sum_{i=1}^{N_0-1}\frac{c_3 i}{\lambda} = \frac{c_3 N_0(N_0-1)}{2\lambda},$$

we obtain the following formula for the mathematical expectation of costs by the joint distribution with event $A_0$

$$E_r(\triangle\gamma_n^{(-)}; A_0) = E_r(\triangle\gamma_n^{(-)} \mid A_0)P(A_0) = \Big(\frac{c_1 N(N+1)}{2\lambda} + c_2(N+N_0) + \frac{c_3 N_0(N_0-1)}{2\lambda} +$$
$$+ c_3 N_0 \int_0^\infty (1 - H_r(y))dy\Big)\int_0^\infty e^{-\lambda y}dH_r(y).$$

Thus, the corresponding mathematical expectation of profit is represented in the following form:

$$E_r(\triangle\gamma_n;A_0) = E_r(\triangle\gamma_n^{(+)};A_0) - E_r(\triangle\gamma_n^{(-)};A_0) = c_0(N+N_0)\int_0^\infty e^{-\lambda y}dH_r(y) - \left[\frac{c_1 N(N+1)}{2\lambda} + \right.$$

$$\left. +c_2(N+N_0) + \frac{c_3 N_0(N_0-1)}{2\lambda} + c_3 N_0\int_0^\infty(1-H_r(y))dy\right]\int_0^\infty e^{-\lambda y}dH_r(y) = \left[c_0(N+N_0) - \right.$$

$$\left. -\frac{c_1 N(N+1)}{2\lambda} - c_2(N+N_0) - \frac{c_3 N_0(N_0-1)}{2\lambda} - c_3 N_0\int_0^\infty(1-H_r(y))dy\right]\int_0^\infty e^{-\lambda y}dH_r(y). \quad (55)$$

8.5.2  Consider the case when s > 0.

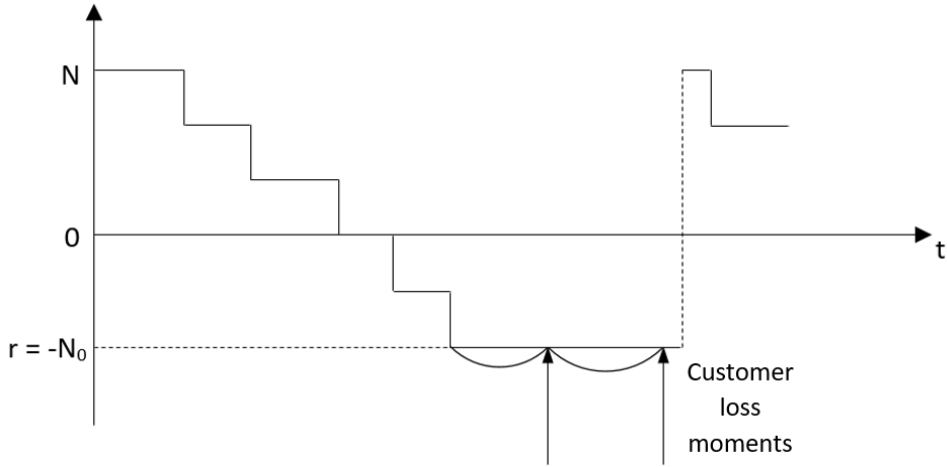

Fig. 24. Graphic illustration for the case 8.5.2

In this case, as in the previous one, income is generated by satisfying the maximum possible number of $N + N_0$ requirements. Consequently,

$$E_r(\triangle\gamma_n^{(+)};A_s) = E_r(\triangle\gamma_n^{(+)}|A_s)P(A_s) = c_0(N+N_0)\int_0^\infty \frac{(\lambda y)^s}{s!}e^{-\lambda y}dH_r(y).$$

Provided that the event $A_s$, $s \geq 1$ is realized, all types of expenses characteristic of the previous case are saved. Analytical expressions for these types of expenses are not changed. The only difference is that there are additional costs associated with the loss of requirements. Since all s requirements received during the delay period will be lost, the amount of these additional costs will be $c_4^{(s)}$. It follows that

$$E_r(\triangle\gamma_n^{(-)}; A_s) = E_r(\triangle\gamma_n^{(-)} | A_s)P(A_s) = \Bigg[\sum_{i=1}^{N}\frac{c_1 i}{\lambda} + \underbrace{c_2(N+N_0)}_{\substack{\text{stock} \\ \text{replenishment} \\ \text{costs}}} + \underbrace{\sum_{i=1}^{N_0-1}\frac{c_3 i}{\lambda}}_{\substack{\text{costs related to the deficit} \\ \text{before the beggining of} \\ \text{the delay period}}} +$$

$$+\underbrace{c_3 N_0 \int_0^\infty (1-H_r(y))dy}_{\substack{\text{costs related to the deficit} \\ \text{after the beggining of} \\ \text{the delay period}}} + \underbrace{c_4^{(s)}}_{\substack{\text{costs due to} \\ \text{loss of clients}}}\Bigg] \int_0^\infty \frac{(\lambda y)^s}{s!}e^{-\lambda y}dH_r(y) = \Bigg[\frac{c_1 N(N+1)}{2\lambda} + c_2(N+N_0) +$$

$$+\frac{c_3 N_0(N_0-1)}{2\lambda} + c_3 N_0 \int_0^\infty (1-H_r(y))dy + c_4^{(s)}\Bigg]\int_0^\infty \frac{(\lambda y)^s}{s!}e^{-\lambda y}dH_r(y).$$

Now we can write out the formula for the mathematical expectation of profit by the joint distribution with the event $A_s$

$$E_r(\triangle\gamma_n; A_s) = E_r(\triangle\gamma_n^{(+)}; A_s) - E_r(\triangle\gamma_n^{(-)}; A_s) = c_0(N+N_0)\int_0^\infty \frac{(\lambda y)^s}{s!}e^{-\lambda y}dH_r(y) - \Bigg[\frac{c_1 N(N+1)}{2\lambda} +$$

$$+c_2(N+N_0) + \frac{c_3 N_0(N_0-1)}{2\lambda} + c_3 N_0 \int_0^\infty (1-H_r(y))dy + c_4^{(s)}\Bigg]\int_0^\infty \frac{(\lambda y)^s}{s!}e^{-\lambda y}dH_r(y) = \Bigg[c_0(N+N_0) -$$

$$-\frac{c_1 N(N+1)}{2\lambda} - c_2(N+N_0) - \frac{c_3 N_0(N_0-1)}{2\lambda} - c_3 N_0 \int_0^\infty (1-H_r(y))dy - c_4^{(s)}\Bigg]\int_0^\infty \frac{(\lambda y)^s}{s!}e^{-\lambda y}dH_r(y). \quad (56)$$

Based on relation (27), taking into account (55), (56), we obtain a representation for the mathematical expectation of profit during the regeneration period, provided that at this period the control parameter takes a fixed value $r = -N_0$

$$E_r(\triangle\gamma_n) = \underbrace{\Bigg[c_0(N+N_0) - \frac{c_1 N(N+1)}{2\lambda} - c_2(N+N_0) - \frac{c_3 N_0(N_0-1)}{2\lambda} -}_{\text{for } s=0}$$

$$\underbrace{-c_3 N_0 \int_0^\infty (1-H_r(y))dy\Bigg]\int_0^\infty e^{-\lambda y}dH_r(y)}_{} + \underbrace{\sum_{i=1}^\infty \Bigg[\Bigg[c_0(N+N_0) - \frac{c_1 N(N+1)}{2\lambda} -}_{\text{for } s>0}$$

$$\underbrace{-c_2(N+N_0) - \frac{c_3 N_0(N_0-1)}{2\lambda} - c_3 N_0 \int_0^\infty (1-H_r(y))dy - c_4^{(i)}\Bigg]\int_0^\infty \frac{(\lambda y)^i}{i!}e^{-\lambda y}dH_r(y)\Bigg]}_{}.$$

At the end of the analysis of this variant, we find an explicit expression for the control efficiency indicator $I_r$, provided that the control parameter takes a fixed value r = -$N_0$. Given that under this condition $E_r(\triangle t_n) = \dfrac{N+N_0}{\lambda} + \int_0^\infty (1-H_r(y))dy$, we obtain

$$I_r = \frac{E_r(\triangle \gamma_n)}{E_r(\triangle t_n)} = \frac{1}{\dfrac{N+N_0}{\lambda} + \int_0^\infty (1-H_r(y))dy} \left\{ \left[ c_0(N+N_0) - \frac{c_1 N(N+1)}{2\lambda} - c_2(N+N_0) - \right.\right.$$

$$\left.- \frac{c_3 N_0(N_0-1)}{2\lambda} - c_3 N_0 \int_0^\infty (1-H_r(y))dy \right] \int_0^\infty e^{-\lambda y} dH_r(y) + \sum_{i=1}^\infty \left[ \left[ c_0(N+N_0) - \frac{c_1 N(N+1)}{2\lambda} - \right.\right.$$

$$\left.\left.- c_2(N+N_0) - \frac{c_3 N_0(N_0-1)}{2\lambda} - c_3 N_0 \int_0^\infty (1-H_r(y))dy - c_4^{(i)} \right] \int_0^\infty \frac{(\lambda y)^i}{i!} e^{-\lambda y} dH_r(y) \right] \right\}. \quad (57)$$

## 9 Solution to the problem of optimal control

The result of the analytical probability study conducted in section 8 is an explicit expression for the function $I_r$ for all permissible values of the parameter $r \in U = \{N, N-1, \ldots, 0, -1, \ldots, -N_0\}$. This explicit expression is defined by formulae (37), (44), (49), (54) and (57). The function $I_r$ coincides with the so-called main function of the linear-fractional integral-type functional $I_\alpha$:

$$I_r = C(r) = \frac{A(r)}{B(r)}, \ r \in U = \{N, N-1, \ldots, 0, -1, \ldots, -N_0\}, \quad (58)$$

In section 6, it was established that optimal control is deterministic and is defined by the point of global extremum (maximum) of the function $I_r$ which is expressed analytically by equality (58). Since the set of permissible controls $U$ is finite, the global maximum of the main function is reached at a certain point $r^* \in U$. Consequently, optimal control exists and coincides with the point of the global maximum $r^*$. The problem of optimal control in the considered stochastic model is reduced to finding point $r^*$, which can be done only numerically.

## 10 Conclusion

Let us make conclusions of the study as a whole. In this work, we have developed a stochastic model of discrete inventory management in the form of a regenerative random process. We have also determined a general approach to solving this problem based on the theorem of unconditional extremum of the linear-fractional integral-type functional depending on discrete probability distribution characterizing the management strategy. Also, in this work, we have obtained explicit analytical expressions for the main function of the linear-fractional integral-type functional that plays the role of the management efficiency indicator. According to the statement of the above-mentioned theorem, the initial management problem can be solved and its solution depends on the global maximum point of the main function defined.